\documentclass[11pt]{article}
\usepackage{definitions}
    
\begin{document}

\title{Nearly Optimal Linear Convergence of Stochastic Primal-Dual Methods for Linear Programming}

\author{Haihao Lu\thanks{The University of Chicago, Booth School of Business (haihao.lu@chicagobooth.edu).} \and Jinwen Yang\thanks{The University of Chicago, Department of Statistics (jinweny@uchicago.edu).}}

\date{December 2023}

\maketitle

\begin{abstract}
  There is a recent interest on first-order methods for linear programming (LP). In this paper, we propose a stochastic algorithm using variance reduction and restarts for solving sharp primal-dual problems such as LP. We show that the proposed stochastic method exhibits a linear convergence rate for solving sharp instances with a high probability.  In addition, we propose an efficient coordinate-based stochastic oracle, which has $\mathcal O(1)$ per iteration cost and improves the complexity of the existing deterministic and stochastic algorithms. Finally, we show that the obtained linear convergence rate is nearly optimal (upto $\log$ terms) for a wide class of stochastic primal-dual methods. Numerical performance verifies the theoretical guarantees of the proposed algorithms.
\end{abstract}

\section{Introduction}

Linear programming (LP), as one of the most fundamental tools in operations research and computer science, has been extensively studied in both academia and industry since 1940s.
The applications of LP span various fields, including pricing and revenue management, transportation, network flow, scheduling, and many others \cite{charnes1954stepping, hanssmann1960linear, bowman1956production,manne1960linear,anderson2000hotel,liu2008choice}. 

The dominant solvers of LP are essentially based on either the simplex method or interior-point method, which are very mature nowadays and can usually provide reliable solutions to LP. However, the success of both methods heavily relies on efficient solving linear systems using factorization, which has the following major drawbacks: (i) the factorization may run out of memory even though the original LP can fit in memory; (ii) it is highly challenging to take advantage of modern computing resources, such as distributed system and GPUs when solving linear systems. These drawbacks make it highly challenging to further scale up LP.

Recently, it was shown that first-order methods (FOMs) with proper enhancements can also identify high-quality solutions to LP problem quickly~\cite{applegate2021practical}, which provides an alternative to the traditional simplex and barrier methods for LP. FOMs utilize only the gradient information of the objective function to update the iterates (in contrast to second-order methods where Hessian is used). The basic operation in FOMs is matrix-vector multiplication, which can avoid the two major drawbacks of traditional methods mentioned above, thus it is suitable for solving larger LP.

In this work, we further push this line of research and study stochastic FOMs for LP. In contrast to deterministic first-order methods~\cite{applegate2021practical}, the basic operation per iteration in stochastic FOMs is a vector operation. Such operation is very cheap in general and is a common practice for modern machine learning applications. As a drawback, the number of iterations of stochastic algorithms to obtain an approximate solution is typically higher than its deterministic counterpart. 

Due to such a fundamental distinction, Nesterov~\cite{nesterov2013gradient} formally defines methods that require matrix factorization (or matrix inversion) as handling \emph{medium-scale} problems, methods that require matrix vector multiplication as handling \emph{large-scale} problems, and methods that require vector operations as handling \emph{huge-scale} problems. Based on this definition, in the context of LP, the simplex method and interior point method are classified as handling medium-scale problems (this definition perhaps belie the practical efficiency of the two methods, but it is a fact that it is challenging to further scale up them as mentioned above), deterministic FOMs~\cite{applegate2021practical,lin2021admm} are classified as handling large-scale problems, and in the paper we instead look at huge-scale problems.

While our motivation is LP, we here consider a more general class of primal-dual problems with the form
\begin{equation}\label{eq:poi}
    \min_{x\in\RR^n}\max_{y\in\RR^m} \mathcal{L}(x,y):= \Phi(x,y) + g_1(x)-g_2(y),
\end{equation}
where $\mathcal L(x,y)$ is convex in $x$ and concave in $y$, $g_1(x)$ is a simple convex function in $x$ and $g_2(y)$ is a simple convex function in $y$.  In particular, the primal-dual formulation of standard form LP is
\begin{equation}\label{eq:lp}
    \min_{x\ge 0}\max_{y\in\RR^m}  y^T Ax + c^Tx-b^Ty \ ,
\end{equation}
and a highly related problem is the unconstrained bilinear problem
\begin{equation}\label{eq:bilinear}
    \min_{x\in\RR^n}\max_{y\in\RR^m}  y^T Ax + c^Tx-b^Ty \ .
\end{equation}
For notational convenience, we define $z=(x,y)$, $F(z)=[\nabla_x\Phi(x,y), -\nabla_y\Phi(x,y)]$ and $g(z)=g_1(x)+g_2(y)$ . We here assume there is an unbiased stochastic oracle $F_{\xi}(z)$ such that $\mathbb E[F_{\xi}(z)]= F(z)$ (see Section \ref{sec:app} for examples on how to construct the stochastic oracles).
To efficiently solve \eqref{eq:lp}, our key ideas are \emph{variance reduction} and \emph{restarts}:

Variance reduction is a successful technique for finite sum stochastic minimization problems~\cite{johnson2013accelerating,defazio2014saga,schmidt2017minimizing}. The basic idea is to reduce the variance in the stochastic gradient estimation by comparing it with the snapshots of the true gradient. Variance reduction can usually improve the convergence property of stochastic minimization algorithms~\cite{johnson2013accelerating,defazio2014saga,schmidt2017minimizing}. Recently, \cite{alacaoglu2021stochastic} extends the variance reduction scheme to EGM for solving monotonic variational inequality (with some ambiguity, we call this algorithm sEGM), which was shown to have $O(1/\epsilon)$ convergence rate. 

Restart is another standard technique for minimization problems, which is used to speed up the convergence of multiple deterministic and stochastic algorithms~\cite{o2015adaptive,freund2018new,yang2018rsg,pokutta2020restarting,zhao2020optimal}. Recently, \cite{applegate2021faster} introduces the sharpness of primal-dual problems, and presents a simple restart scheme to accelerate the linear convergence rate of primal-dual algorithms for solving sharp problems.

This paper extends the restarted algorithm in \cite{applegate2021faster} to stochastic algorithms (in particular sEGM~\cite{alacaoglu2021stochastic}) for solving sharp primal-dual problems, such as LP. We further show that the obtained linear convergence rate is nearly optimal (upto $\log$ terms) for a wild class of stochastic algorithms. It turns out that while LP is sharp on any bounded region~\cite{applegate2021faster}, LP is not globally sharp (see Appendix~\ref{sec:sharp-subreg} and Appendix~\ref{sec:lp_not_global_sharp} for a counter example). A fundamental difficulty of restarted stochastic algorithm is that, unlike deterministic algorithms, the iterates may escape from any bounded region, thus there is no guarantee that local sharpness can be helpful for the convergence of stochastic algorithms. We overcome this issue by presenting a high-probability argument and show the linear convergence of the proposed restarted algorithms with high probability. 

The performance of randomized algorithms depends on the realization of stochasticity. The traditional convergence analysis of these algorithms usually measures the expected performance, namely, running the algorithm multiple times and looking at the average performance. In contrast, we present high probability results, namely, running the algorithm multiple times, and study the performance of all trajectories with a certain probability. In particular, when choosing the probability as $1/2$, we obtain the convergence guarantee of the median trajectory, which can be viewed as an alternative to the expected performance studied in the related literature.

\begin{table}
\centering
\small{
\begin{threeparttable}
\begin{tabular}{|c|c|c|c|}
     \hline
     \thead{Algorithm} & \thead{\makecell{Per-iteration\\ Cost}} & \thead{Number of Iteration to\\ find an $\epsilon$ solution} & \thead{Total Cost\footnotemark} \\ \hline
     \makecell{Deterministic \\  \cite{korpelevich1976extragradient}} & $\mathcal{O}\pran{\text{nnz}(A)}$ & $\mathcal{O}\pran{\frac{\Vert A \Vert_2}{\epsilon}}$ & $\mathcal{O}\pran{\text{nnz}(A)\frac{\Vert A \Vert_2}{\epsilon}}$ \\ \hline
     \makecell{Deterministic \\ Restart  \cite{applegate2021faster}} & $\mathcal{O}\pran{\text{nnz}(A)}$ & $\mathcal{O}\pran{\kappa_2 \log{\frac{1}{\epsilon}}}$ & $\mathcal{O}\pran{\text{nnz}(A)\kappa_2 \log{\frac{1}{\epsilon}}}$ \\ \hline
     \makecell{SPDHG \\ \cite{chambolle2018stochastic}} & $\mathcal{O}(m+n)$ & ${\mathcal{O}}\pran{\frac{\sum_i \Vert A_{i\cdot} \Vert_2}{\epsilon}}$ & ${\mathcal{O}}\pran{\text{nnz}(A)+(m+n)\frac{\sum_i \Vert A_{i\cdot} \Vert_2}{\epsilon}}$ \\ \hline
     \makecell{SPDHG\footnotemark \\ \cite{alacaoglu2019convergence}} & $\mathcal {O}(m+n)$ & ${\mathcal O}\pran{\kappa_F^2 \log \frac 1 \epsilon}$ & ${\mathcal O}\pran{(m+n)\kappa_F^2 \log \frac 1 \epsilon}$ \\ \hline
     \makecell{Conceptual \\ Proximal \cite{carmon2019variance}} & $\mathcal{O}(m+n)$ & $\mathcal{O}\pran{\sqrt{\frac{\text{nnz}(A)}{m+n}}\frac{\Vert A \Vert_F}{\epsilon}}$ & \makecell{$\mathcal{O}\pran{\text{nnz}(A)+\frac{\sqrt{\text{nnz}(A)(m+n)}\Vert A \Vert_F}{\epsilon}}$} \\  \hline
     \makecell{sEGM \\ \cite{alacaoglu2021stochastic}} & $\mathcal{O}(m+n)$ & $\mathcal{O}\pran{\sqrt{\frac{\text{nnz}(A)}{m+n}}\frac{\Vert A \Vert_F}{\epsilon}}$ & \makecell{$\mathcal{O}\pran{\text{nnz}(A)+\frac{\sqrt{\text{nnz}(A)(m+n)}\Vert A \Vert_F}{\epsilon}}$} \\  \hline
     \makecell{RsEGM\\ Oracle II \\ {(\bf This Paper)}} & $\mathcal{O}(m+n)$ & $\widetilde{\mathcal{O}}\pran{\sqrt{\frac{\text{nnz}(A)}{m+n}}\kappa_F\text{polylog}(\frac 1 \epsilon)}$ &  \makecell{$\widetilde{\mathcal O}\pran{\text{nnz}(A) \log{\frac{1}{\epsilon}}+\sqrt{\text{nnz}(A)(m+n)}\kappa_F\text{polylog}(\frac 1 \epsilon)}$ } \\ \hline
     \makecell{RsEGM\\ Oracle IV \\ {(\bf This Paper)}} & $\mathcal{O}(1)$ & $\widetilde{\mathcal{O}}\pran{\sqrt{{\text{nnz}(A)}}\kappa_F\text{polylog}(\frac 1 \epsilon)}$ &  \makecell{$\widetilde{\mathcal O}\pran{\text{nnz}(A) \log{\frac{1}{\epsilon}}+\sqrt{\text{nnz}(A)}\kappa_F\text{polylog}(\frac 1 \epsilon)}$ \\  } \\ \hline
\end{tabular}
\caption{Comparison on Unconstrained Bilinear Problem \eqref{eq:bilinear}, where $\kappa_2=\frac{\Vert A \Vert_2}{\alpha}$, $\kappa_F=\frac{\Vert A \Vert_F}{\alpha}$.}
\label{tab:bilinear}
\end{threeparttable}}
\end{table}

\footnotetext[1]{
{The total cost is sometimes more than the product of the per-iteration cost and the number of iterations due to the snapshot step in variance reduction.}}
\footnotetext[2]{
The complexity of SPDHG~\cite{alacaoglu2021stochastic} involves complicated terms. In the table, we present a lower bound on the number of iterations and the total cost, i.e., they need at least this number of iterations (or total cost) to find an $\epsilon$-accuracy solution based on their analysis (see Appendix \ref{sec:compare-spdhg} for more details).}

\begin{table}
\centering
{\small
\begin{tabular}{ |c|c|c|c| }
     \hline
     \thead{Algorithm} & \thead{\makecell{Per-iteration\\ Cost}} & \thead{Number of Iteration} & \thead{Total Cost} \\ \hline
     \makecell{Deterministic \\  \cite{korpelevich1976extragradient}} & $\mathcal{O}(\text{nnz}(A))$ & $\mathcal{O}\pran{\frac{\Vert A \Vert_2}{\epsilon}}$ & $\mathcal{O}\pran{\text{nnz}(A)\frac{\Vert A \Vert_2}{\epsilon}}$ \\ \hline
     \makecell{Deterministic\\ Restart \cite{applegate2021faster}} & $\mathcal{O}(\text{nnz}(A))$ & $\mathcal{O}\pran{\kappa_2\log{\frac{1}{\epsilon}}}$ & $\mathcal{O}\pran{\text{nnz}(A)\kappa_2\log{\frac{1}{\epsilon}}}$ \\ \hline
     \makecell{sEGM \\ \cite{alacaoglu2021stochastic}} & $\mathcal{O}(m+n)$ & $\mathcal{O}\pran{\sqrt{\frac{\text{nnz}(A)}{m+n}}\frac{\Vert A \Vert_F}{\epsilon}}$ & \makecell{$\mathcal{O}\pran{\text{nnz}(A)+\frac{\sqrt{\text{nnz}(A)(m+n)}\Vert A \Vert_F}{\epsilon}}$} \\  \hline
     \makecell{RsEGM \\ Oracle II \\ {(\bf This Paper)}} & $\mathcal{O}(m+n)$ & $\widetilde{\mathcal{O}}\pran{\sqrt{\frac{\text{nnz}(A)}{m+n}}\kappa_F\text{polylog}(\frac 1 \epsilon)}$ &  \makecell{$\widetilde{\mathcal O}\pran{\text{nnz}(A)\log{\frac{1}{\epsilon}}+\sqrt{\text{nnz}(A)(m+n)}\kappa_F\text{polylog}(\frac 1 \epsilon)}$} \\ \hline
     \makecell{RsEGM \\ Oracle IV \\ {(\bf This Paper)}} & $\mathcal{O}(1)$ & $\widetilde{\mathcal{O}}\pran{\sqrt{{\text{nnz}(A)}}\kappa_F\text{polylog}(\frac 1 \epsilon)}$ &  \makecell{{ $\widetilde{\mathcal O}\pran{\text{nnz}(A)\log{\frac{1}{\epsilon}}+\sqrt{\text{nnz}(A)}\kappa_F\text{polylog}(\frac 1 \epsilon)}$} }  \\ \hline
\end{tabular}}
\caption{Comparison on standard LP~\eqref{eq:lp}, where $\kappa_2=\frac{\Vert A \Vert_2}{1/[H(1+\Vert z^{0,0} \Vert+R_0)]}$, $\kappa_F=\frac{\Vert A \Vert_F}{1/[H(1+\Vert z^{0,0} \Vert+R_0)]}$, and $H$ is the Hoffman constant of the KKT system of the LP (see Example \ref{thm:lem-LP-sharp} for details).}
\label{tab:lp}
\end{table}

Table \ref{tab:bilinear} and Table \ref{tab:lp} present the per iteration cost, complexity of the algorithm, and the total flop counts for our proposed algorithm, Restarted sEGM (RsEGM), and compare it with multiple deterministic and stochastic algorithms for solving unconstrained bilinear problem~\eqref{eq:bilinear} and LP~\eqref{eq:lp}, respectively. For unconstrained bilinear problems~\eqref{eq:bilinear}, deterministic algorithms require to compute the gradient of the objective, thus the per iteration cost is $O(\text{nnz}(A))$. Standard stochastic algorithms are usually based on row/column sampling, and the iteration cost is $O(m+n)$. We also present a stochastic coordinate scheme (oracle IV) where the per-iteration cost is $O(1)$. While stochastic algorithms have low per-iteration cost, it usually requires more iterations to identify an $\epsilon$-close solution compared to its deterministic counterparts. As we see in Table \ref{tab:bilinear}, compared with the optimal deterministic algorithms~\cite{applegate2021faster}, the total flop count of RsEGM with stochastic Oracle II is better when the matrix $A$ is dense and of low rank. With the coordinate gradient estimator Oracle IV, the total cost of RsEGM is even lower and improves optimal deterministic algorithms by at least a factor of $\sqrt n$ when A is dense. For standard LP, as seen in Table \ref{tab:lp}, stochastic algorithms require more iterations to achieve $\epsilon$-accuracy than the unconstrained bilinear problem due to the existence of inequality constrains. Similar to the unconstrained bilinear setting, when matrix $A$ is low-rank and dense, the total flop cost of RsEGM improves the optimal deterministic algorithm. On the other hand, most of the previous works on stochastic algorithms have sublinear rate. The only exception is~\cite{alacaoglu2021stochastic}, where the authors show the linear convergence of SPDHG for solving problems satisfying global metric sub-regularity. Indeed, unconstrained bilinear problems satisfy the global metric sub-regularity, while LP does not satisfy it globally. The complexity of SPDHG involves more complicated notations and we present a more detailed comparison in Appendix~\ref{sec:lp_not_global_sharp} and Appendix~\ref{sec:compare-spdhg}, but our proposed algorithms are at least better than SPDHG by a factor of condition number $\kappa_F$.









\subsection{Summary of Contributions}

The contributions of the paper can be summarized as follow:

    (i) We propose a restarted stochastic extragradient method with variance reduction for solving sharp primal-dual problems. We show that the proposed algorithm exhibit linear convergence rate with high probability, in particular,
    \begin{itemize}
        \item For linear programming and unconstrained bilinear problems, our restarted scheme improves the complexity of existing linear convergence of stochastic algorithms~\cite{alacaoglu2021stochastic} by a factor of the condition number. The improvement comes from restarts.
        \item To the best of our knowledge, this is the first stochastic algorithm with linear rate for the general standard-form LP problems~\eqref{eq:lp}. To prove this result, we introduce a high-probability analysis to upper-bound distance between the iterates and the optimal solution set.
    \end{itemize}
     
    (ii) We present the complexity lower bound of a class of stochastic first-order methods for solving sharp primal-dual problems, which matches the upper bound we obtained for RsEGM (upto $\log$ terms). This showcases RsEGM achieves a nearly optimal linear convergence rate for sharp primal-dual problems.





\subsection{Assumptions}
Throughout the paper, we have two assumptions on the problem and the stochastic oracle. The first one is on the primal-dual problem: 
\begin{ass}\label{ass:problem}
The problem \eqref{eq:poi} satisfies:\\
    (i) $\mathcal L(x,y)$ is convex in $x$ and concave in $y$.\\
    (ii) $g: \mathcal Z \rightarrow \mathbb R \cup {+\infty}$ is proper convex lower semi-continuous.\\
    (iii) The stationary solution set $\mathcal Z^*=\{z|0\in F(z)+\partial g(z)\} \neq \varnothing$.
\end{ass}

As a stochastic first-order method, we assume there exists a stochastic gradient oracle:
\begin{ass}\label{ass:oracal}
We assume there exists a stochastic oracle $F_\xi:\RR^{m+n}\rightarrow \RR^{m+n}$ such that

(i) it is unbiased: $\mathbb E[F_{\xi}(z)]=F(z)$;

(ii) it is $L$-Lipschitz (in expectation): $\mathbb E [\Vert F_{\xi}(u)-F_{\xi}(v) \Vert ^2]\leq L^2\Vert u-v \Vert ^2$.

\end{ass}

\subsection{Related Literature}
~\\ \textbf{Convex-concave primal-dual problems.}
There has been a long history of convex-concave primal-dual problems, and many of the early works study a more general problem, monotone variational inequalities. Rockafellar proposed a proximal point method (PPM) \cite{rockafellar1976monotone} for solving monotone variational inequalities. Around the same time, Korpelevich proposed the extragradient method (EGM)  \cite{korpelevich1976extragradient} for convex-concave primal-dual problems. After that, there have been numerous results on the convergence analysis of these methods. In particular, Tseng \cite{tseng1995linear} shows that PPM and EGM have linear convergence for strongly-convex-strongly-concave primal-dual problems or for unconstrained bilinear problems. Nemirovski proposes Mirror Prox algorithm in the seminal work \cite{nemirovski2004prox}, which is a more general form of EGM, and shows that EGM has $\mathcal O(\frac 1 \epsilon)$ sublinear convergence rate for solving general convex-concave primal-dual problems over a bounded and compact set.  \cite{nemirovski2004prox} also build up the connection between EGM and PPM: EGM is an approximation to PPM.

Another line of research is to study a special case of \eqref{eq:poi} where $\Phi(x,y)=y^TAx$ is a bilinear term. Two well-known algorithms are Douglas-Rachford splitting \cite{douglas1956numerical,eckstein1992douglas} (Alternating Direction Method of Multiplier (ADMM) as a special case) and Primal-dual Hybrid Gradient Method (PDHG) \cite{chambolle2011first}.

Very recently, there is a renewed interest on primal-dual methods, motivated by machine learning applications. For bilinear problems $\mathcal L(x,y)=y^TAx$ with full rank matrix $A$,  \cite{daskalakis2018training} shows that the Optimistic Gradient Descent Ascent (OGDA) converges linearly and later on \cite{mokhtari2020unified} shows that OGDA , EGM and PPM all enjoy a linear convergence rate. \cite{mokhtari2020unified} also presents an interesting observation that OGDA approximates PPM on bilinear problems. Lu \cite{lu2020s} analyzes the dynamics of unconstrained primal-dual algorithms under an ODE framework and yields tight conditions under which different algorithms exhibit linear convergence. However, an important caveat is that not all linear convergence rates are equal. \cite{applegate2021faster} shows that a simple restarted variant of these algorithms can improve the dependence of complexity on condition number in their linear convergence rate, as well as the empirical performance of the algorithms. 

\textbf{Linear programming.}
Linear programming is a fundamental tool in operations research. Two dominating methods to solve LP problems are simplex method \cite{dantzig1998linear} and interior-point method \cite{karmarkar1984new}, and the commercial LP solvers based on these methods can provide reliable solutions even for fairly large instances. While the two methods are quite different, both require solving linear systems using factorization. As a result, it becomes very challenging to further scale up these two methods, in particular to take advantage of distributed computing. Recently, there is a recent trend on developing first-order methods for LP only utilizing matrix-vector multiplication \cite{pmlr-v119-basu20a,ipm25gondzio2012,lin2021admm,applegate2021practical,applegate2021infeasibility}. In general, these methods are easy to be parallelized and do not need to store the factorization in memory. 

The traditional results of first-order methods for LP usually have sublinear rate, due to the lack of strong convexity, which prevents them to identify high-accuracy solutions. To deal with this issue, \cite{eckstein1990alternating} presents a variant of ADMM and shows the linear convergence of the proposed method for LP. More recently \cite{lewis2018partial,liang2018local} show that many primal-dual algorithms under a mild non-degeneracy condition have eventual linear convergence, but it may take a long time before reaching the linear convergence regime. \cite{applegate2021faster} propose a restarted scheme for LP in the primal-dual formulation. They introduce a sharpness condition for primal-dual problems based on the normalized duality gap and show that the primal-dual formulation of LP is sharp on any bounded region. Then they provide restarted schemes for sharp primal-dual problems and show that their proposed algorithms have the optimal linear convergence rate (in a class of deterministic first-order methods) when solving sharp problems. A concurrent work \cite{song2022coordinate} studies stochastic algorithms for generalized linear programming, but the focus is on sublinear rate, while our focus is on the linear rate.

\textbf{Sharpness conditions and restart schemes.} The concept of sharpness was first proposed by Polyak \cite{polyak1979sharp} on minimization problem. Recently, there is a trend of work on developing first-order method with faster convergence rates using sharpness. For example, \cite{yang2018rsg} shows linear convergence of restarted subgradient descent on sharp non-smooth functions and there are other works on sharp non-convex minimization \cite{davis2018subgradient}. Sharpness can also be viewed as a certain error bound condition \cite{roulet2020sharpness}. Recently \cite{applegate2021faster} introduces sharpness condition for primal-dual problems. A highly related concept is the metric subregularity for variational inequalities, which is a weaker condition than the sharpness condition proposed in \cite{applegate2021faster} (see Appendix \ref{sec:sharp-subreg} for a discussion). Under such conditions, \cite{alacaoglu2019convergence,fercoq2021quadratic} present the linear convergence for stochastic PDHG and deterministic PDHG, respectively.


Restarting is a powerful technique in optimization. It can improve the practical and theoretical convergence of a base algorithm without modification to the base algorithm \cite{pokutta2020restarting}. Recently, there have been extensive works on this technique in smooth convex optimization \cite{o2015adaptive,roulet2020sharpness}, non-smooth convex optimization \cite{freund2018new,yang2018rsg} and stochastic convex optimization \cite{johnson2013accelerating,lin2015universal,tang2018rest}. For sharp primal-dual problems, \cite{applegate2021faster} propose fixed-frequency and adaptive restart on a large class of base primal-dual algorithms including PDHG, ADMM and EGM.

\textbf{Variance reduction and primal-dual problem.}
Variance reduction technique \cite{defazio2014saga,johnson2013accelerating} is developed to improve the convergence rate for stochastic algorithms upon pure SGD for minimization problems. There are extensive works on variants of stochastic variance reduction for minimization problems under various settings (see \cite{gower2020variance} for a recent overview).

Compared with the extensive works on minimization problem, the research of variance-reduced methods on primal-dual problems is fairly limited. \cite{palaniappan2016stochastic} studies stochastic forward-backward algorithm with variance reduction for primal-dual problems and, more generally, monotone inclusions. Under strong monotonicity, they prove a linear convergence rate and improve the complexity of deterministic methods for bilinear problems. \cite{carmon2019variance} proposes a randomized variant of Mirror Prox algorithm. They focus on matrix games and improve complexity over deterministic methods in several settings. However, beyond matrix games, their method requires extra assumptions such as bounded domain and involves a three-loop algorithm. More recently, \cite{alacaoglu2021stochastic} proposes a stochastic extragradient method with variance reduction for solving variational inequalities. Under Euclidean setting, their method is based on a loopless variant of variance-reduced method \cite{hofmann2015variance,kovalev2020don}. Their algorithm offers a similar convergence guarantee as \cite{carmon2019variance} but does not require assumptions such as bounded region. 

{After the first version of this paper, there have been several recent works on stochastic primal-dual methods that can be applied to LP. For example, \cite{huang2022accelerated,nan2023extragradient} propose variants of variance-reduced extra-gradient methods for finite-sum variational inequality. In particular, \cite{huang2022accelerated} derives the sublinear rate of the proposed method under the finite-sum convex-concave setting, while \cite{nan2023extragradient} derives the linear convergence rate of Algorithm \ref{alg:segm} under the projection-type error-bound condition. \cite{alacaoglu2022complexity} studies complexity bounds for the primal-dual algorithm with random extrapolation and coordinate descent and obtains sub-linear rate for general convex-concave problems with bilinear coupling. }

\subsection{Notations}

Let $\log(\cdot)$ and $\exp(\cdot)$ refer to the natural log and exponential function, respectively. Let $\sigma_{min}(A)$ and $\sigma_{min}^+(A)$ denote the minimum singular value and minimum nonzero singular value of a matrix $A$. Let $\Vert z \Vert_2=\sqrt{\sum_j z_j^2}$ be the $l_2$ norm of vector $z$. Let $\Vert A \Vert_2$ and $\Vert A \Vert_F$ denote the 2-norm and Frobenius norm of a matrix $A$. Let $e_i$ be the $i$-th standard unit vector. Let $f(x)=\mathcal{O}(g(x))$ denote for sufficiently large $x$, there exists constant $C$ such that $f(x)\leq Cg(x)$ and $f(x)=\Omega(g(x))$ denote for sufficiently large $x$, there exists constant $c$ such that $f(x)\geq cg(x)$. $f(x)=\Theta(g(x))$ if $f(x)=\mathcal{O}(g(x))$ and $f(x)=\Omega(g(x))$. Denote $\widetilde{\mathcal{O}}(g(x))=\mathcal{O}\bigg(g(x)\text{polylog}(g(x))\bigg)$, where $\text{polylog}(g(x))$ refers to a polynomial in $\log(g(x))$. For set $\mathcal Z$, denote $W_r(z)=\left\{\hat z\in \mathcal Z \big| \Vert z-\hat z \Vert\leq r\right\}$ the ball centered at $z$ with radius $r\in (0,\infty)$ intersected with the set $\mathcal Z$. $\iota_S$ is the indicator function of convex set $S$.

\section{Preliminaries}

In this section, we present two results in the recent literature (i) the sharpness condition for primal-dual problems based on the normalized duality gap introduced in \cite{applegate2021faster}, and (ii) the stochastic EGM with variance reduction that was introduced in \cite{alacaoglu2021stochastic}. We will utilize these results to develop our main theory in later sections.

\subsection{Sharpness Condition for Primal-Dual Problems}\label{sec:sharpness}
Sharpness is a central property of the objective in minimization problems that can speed up the convergence of optimization methods. 
Recently, \cite{applegate2021faster} introduces a new sharpness condition of primal-dual problems based on a normalized duality gap. They show that linear programming (among other examples) in the primal-dual formulation is sharp. Furthermore, \cite{applegate2021faster} proposes a restarted algorithm that can accelerate the convergence of primal-dual algorithms such as EGM, PDHG and ADMM, and achieve the optimal linear convergence rate on sharp primal-dual problems.

More formally, the sharpness of a primal-dual problem is defined as:
\begin{mydef}
    We say a primal-dual problem is $\alpha$-sharp on the set $S\subset \mathcal Z$ if $\rho_r(z)$ is $\alpha$-sharp on $S$ for all $r$, i.e. it holds for all $z\in S$ that 
    \begin{equation*}
        \alpha\text{dist}(z,\mathcal Z^*)\leq \rho_r(z)=max_{\hat z\in W_r(z)}\frac {\mathcal L(x,\hat y)-\mathcal L(\hat x,y)}{r}\ ,
    \end{equation*}
    where $W_r(z)=\left\{\hat z\in \mathcal Z \big| \Vert z-\hat z \Vert\leq r\right\}$ is a ball centered at $z$ with radius $r$ intersected with the set $\mathcal Z$.
\end{mydef}

In particular, the following examples are shown to be sharp instances~\cite{applegate2021faster}:
\begin{exam}\label{exam:bd}
    Consider a primal-dual problem with a bounded feasible region, i.e., $g_1$ and $g_2$ encode a bounded feasible region of $x$ and $y$, respectively. Suppose $P(x)=\max_{y}\mathcal L(x,y)$ is $\alpha_1$-sharp in $x$ and $D(y)=\max_{x\in \mathcal X}\mathcal L(x,y)$ is $\alpha_2$-sharp in $y$. Then problem \eqref{eq:poi} is $\alpha$-sharp in the feasible region, where $\alpha=\frac{\min{\{\alpha_1,\alpha_2\}}}{R}$ and $R$ is the diameter of the region.
\end{exam} 

\begin{exam}\label{exam:bilinear}
    Consider an unconstrained bilinear problem with $\mathcal L(x,y)=y^TAx+c^Tx-b^Ty$, then problem \eqref{eq:poi} is $\alpha$-sharp in $\mathbb R^{m+n}$, where $\alpha=\sigma^+_{min}(A)$.
\end{exam}

\begin{exam}\label{thm:lem-LP-sharp}
    Consider the primal-dual form of standard LP: $$\min_{x\in \mathbb R^n}\max_{y\in \mathbb R^m} \mathcal L(x,y)=y^TAx+c^Tx+\iota_{x\geq 0}-b^Ty\ .$$ Then for any $z$, the normalized duality gap $\rho_r(z)=\max_{\hat z\in W_r(z)}\frac {\mathcal L(x,\hat y)-\mathcal L(\hat x,y)}{r}$ satisfies $$\rho_r(z)\geq \frac{1}{H\sqrt{1+4\Vert z \Vert ^2}}\text{dist}(z,\mathcal Z^*)\ ,$$ where $H$ is the Hoffman constant of the KKT system of LP \cite{applegate2021faster,hoffman1952approximate}.
\end{exam}

\subsection{Stochastic EGM with Variance Reduction}

In this subsection, we present the stochastic EGM with variance reduction introduced in \cite{alacaoglu2021stochastic}, and restate their two major convergence results. 


Algorithm \ref{alg:segm} presents the stochastic EGM with variance reduction (sEGM). The algorithm keeps track of two sequences $\{z_k\}$ and $\{w_k\}$. $w_k$ can be viewed as a snapshot, which is updated rarely and we evaluate the full gradient $F(w_k)$. Then we calculate $\bz_k$ as a weighted average of $z_k$ and $w_k$. Similar to deterministic EGM, the algorithm computes an intermediate solution $z_{k+1/2}$. We then use the full gradient of the snapshot $F(w_k)$ to compute $z_{k+1/2}$, and use the variance reduced stochastic gradient $F(w_k)+F_{\xi_k}(z_{k+1/2})-F_{\xi_k}(w_k)$ to compute $z_{k+1}$.

\begin{algorithm}
    \renewcommand{\algorithmicrequire}{\textbf{Input:}}
    \renewcommand{\algorithmicensure}{\textbf{Output:}}
    \caption{Stochastic EGM with variance reduction: $z^K=sEGM(p,Q,\tau,z^0,K)$}
    \label{alg:segm}
    \begin{algorithmic}[1]
        \REQUIRE Probability $p\in(0,1)$, sample oracle $Q$, step size $\tau$, number of iteration $K$. Let $w^0=z^0$.
        \ENSURE $\widetilde z^{K}=\frac 1 K \sum_{l=0}^{k-1}z^{l+1/2}$.
        \FOR{$k=0,1,...K-1$}
        \STATE $\bar z^k=(1-p) z^k + p w^k$
        \STATE $z^{k+1/2}=\text{prox}_{\tau g}(\bar z^k-\tau F(w^k))$
        \STATE Draw an index $\xi_k$ according to Q
        \STATE $z^{k+1}=\text{prox}_{\tau g}(\bar z^k-\tau[F(w^k)+F_{\xi_k}(z^{k+1/2})-F_{\xi_k}(w^k)])$
        \STATE $w^{k+1}= \begin{cases} z^{k+1}, & \text{with prob.}\; p\\ w^{k}, & \text{with prob.}\; 1-p \end{cases} $
        \ENDFOR
    \end{algorithmic}
\end{algorithm}

We here restate the descent lemma of Algorithm \ref{alg:segm} (\cite[Lemma 2.2]{alacaoglu2021stochastic}), which is the key component in the sublinear convergence rate proof of sEGM:
\begin{lem}\label{lem:lem1}
    Let $\phi_k(z)=(1-p) \Vert z^k-z\Vert ^2+\Vert w^k-z \Vert^2$. Let $\tau=\frac{\sqrt{p}}{2L}$. Then for $(z^k)$ generated by Algorithm \ref{alg:segm} and any $z_{*}\in \mathcal Z^{*}$, it holds that 
    \begin{equation*}
        \mathbb{E}_k[\phi_{k+1}(z_*)]\leq \phi_k(z_*)-\frac 12\pran{p\Vert z^{k+1/2}-w_k\Vert ^2+\mathbb{E}_k[\Vert z^{k+1/2}-z^{k+1}\Vert^2]}\ .
    \end{equation*}
    Moreover, it holds that 
    \begin{equation*}
        \sum_{k=0}^{\infty}\pran{p\mathbb{E}[\Vert z^{k+1/2}-w^{k}\Vert^2]+\mathbb{E}[\Vert z^{k+1/2}-z^{k+1}\Vert^2]}\leq 2\phi_0(z_*)=2(2-p)\Vert z^0-z_* \Vert^2\ .
    \end{equation*}
\end{lem}


\vspace{0.2cm}
Next, we restate the sublinear convergence of Algorithm \ref{alg:segm} in terms of the duality gap. The theorem and its proof can be obtained by a simple modification of \cite[Theorem 2.5]{alacaoglu2021stochastic}.
\begin{thm}\label{thm:thm1}
    Let Assumptions hold, $p\in (0,1)$ and $\tau=\frac{\sqrt{p}}{2L}$. 
    Let $\Theta_{k+1/2}(z)=\langle F(z^{k+1/2}),z^{k+1/2}-z \rangle +g(z^{k+1/2})-g(z)$. Then for any compact set $\mathcal C$,
    \begin{equation*}
        2\tau\mathbb{E}\left[\max_{z\in \mathcal C}\sum_{l=0}^{k-1} \Theta_{l+1/2}(z)\right]\leq \frac 72 \max_{z\in \mathcal C}\Vert z^0-z\Vert^2+14\Vert z^0-z_*\Vert^2\ .
    \end{equation*}
\end{thm}

\section{Restarted Stochastic EGM with Variance Reduction}\label{sec:convergence}
In this section, we first present our main algorithm, the nested-loop restarted stochastic EGM algorithm (RsEGM, Algorithm \ref{alg:rsegm}), in section \ref{sec:algorithm}. Then we show that with high probability, the restarted algorithm exhibits linear convergence to an optimal solution for sharp primal-dual problems both with and without bounded region. This linear convergence result accelerates the linear convergence of stochastic algorithms without restart.


\subsection{Algorithm}\label{sec:algorithm}
Algorithm \ref{alg:rsegm} presents the nested-loop restarted stochastic EGM algorithm (RsEGM). We initialize the algorithm with probability parameter $p\in(0,1]$, a stochastic sample oracle $Q$, step size $\tau$ for sEGM, length of inner loop $K$ and outer iteration $T$. In each outer loop, we run sEGM for $K$ iterations and restart the next outer loop using the output of sEGM from the previous outer loop. 



\begin{algorithm} 
    \renewcommand{\algorithmicrequire}{\textbf{Input:}}
    \renewcommand{\algorithmicensure}{\textbf{Output:}}
    \caption{Restarted stochastic EGM: $z^{T,0}=RsEGM(p,Q,\tau,z^{0,0},T,K)$}
    \label{alg:rsegm}
    \begin{algorithmic}[1]
        \REQUIRE Probability $p\in(0,1]$, probability distribution $Q$, step size $\tau$, number of restart $T$, number of iteration for sEGM $K$.
        \ENSURE $z^{T,0}$.
        \FOR{$t=0,1,...,T-1$}
        \STATE Call subroutine sEGM to obtain $z^{t+1,0}=sEGM(p,Q,\tau,z^{t,0}, K)$
        \ENDFOR
    \end{algorithmic}
\end{algorithm}

\subsection{Convergence Guarantees for Problems with Global Sharpness}\label{sec:global-sharp}

Theorem \ref{thm:thm-global} presents the high probability linear convergence rate of Algorithm \ref{alg:rsegm} for solving primal-dual problem \eqref{eq:poi} with global sharpness conditions. Recall that the global sharpness holds for unconstrained bilinear problems and bilinear problems with bounded region (an economic example is two-player matrix game), this implies the linear convergence of Algorithm \ref{alg:rsegm} for these two cases. 



\begin{thm}\label{thm:thm-global}
    Consider RsEGM (Algorithm \ref{alg:rsegm}) for solving a sharp primal-dual problem~\eqref{eq:poi}. Suppose Assumption \ref{ass:problem} and Assumption \ref{ass:oracal} hold. Let $\alpha$ be the global sharpness constant of the primal-dual problem and $L$ be the Lipschitz parameter of the stochastic oracle. Denote the initial distance to optimal set as $R_0=\text{dist}(z^{0,0},\mathcal Z^*)$, and choose step-size $\tau=\frac{\sqrt{p}}{2L}$.  For any $\delta\in (0,1)$ and $\epsilon>0$, if the outer loop count $T$ and the inner loop count $K$ of RsEGM satisfy $$T\ge \max{\left\{ \log_{2}\frac{R_0}{\epsilon}, 0 \right\}},\; K\geq \widetilde{\mathcal{O}}\pran{\frac{1}{\sqrt p}\frac{L}{\alpha}\frac{1}{\delta^2}T^2}\ ,$$ then it holds with probability at least $1-\delta$ that $$\text{dist} (z^{T,0},\mathcal Z^*)\leq \epsilon\ .$$
    Furthermore, the total number of iteration to get $\epsilon$ accuracy with probability at least $1-\delta$ is of order $$\widetilde{\mathcal{O}}\left(\frac{1}{\sqrt p}\frac{L}{\alpha}\frac{1}{\delta^2}\pran{\log\frac{R_0}{\epsilon}}^3\right)\ .$$ 
\end{thm}

\begin{rem}\label{rem:median}
    The performance of stochastic algorithms depends on the realization of stochasticity. The traditional convergence analysis of these algorithms usually measures the expected performance, namely, running the algorithm multiple times and looking at the average performance. In contrast, when choosing $\delta=\frac{1}{2}$, Theorem \ref{thm:thm-global} provides the complexity of the median performance of the algorithm.
\end{rem}

\begin{rem}\label{rem:finite-sum}
    In this remark, we consider the finite-sum model where $\mathcal L(x,y)=\frac{1}{N}\sum_{i=1}^N \mathcal L_i(x,y)$. Suppose the stochastic oracle computes $F_{\xi}(z)=[\nabla_x \mathcal L_i(x,y), -\nabla_y \mathcal L_i(x,y)]$ for $i$ that is uniformly random in set $[1,...,N]$. Then evaluating $F_{\xi}(z)$ needs one stochastic oracle and evaluating the true gradient $F(z)$ need $N$ stochastic oracles. With a careful calculation, we obtain that one outer iteration of RsEGM requires $\pran{N+2}+\widetilde{\mathcal{O}}\pran{\pran{Np+2}\frac{1}{\sqrt p}\frac{L}{\alpha}\frac{1}{\delta^2}\pran{\log{\frac{R_0}{\epsilon}}}^2}$ stochastic oracle calls. Thus, the average total oracle cost of RsEGM to reach $\epsilon$-accuracy with probability at least $1-\delta$ becomes 
    \begin{equation}\label{eq:bd-oracle-complexity}
        \widetilde{\mathcal{O}}\pran{\pran{N+2}\log{\frac{R_0}{\epsilon}}+\pran{Np+2}\frac{1}{\sqrt p}\frac{L}{\alpha}\frac{1}{\delta^2}\pran{\log{\frac{R_0}{\epsilon}}}^3}\ .
    \end{equation} Choosing $p=\Theta{\pran{\frac{1}{N}}}$ and $\delta=0.5$, the number of stochastic oracle calls for the median trajectory to achieve $\epsilon$ accuracy becomes $$\widetilde{\mathcal{O}}\pran{\pran{N+2}\log{\frac{R_0}{\epsilon}}+\sqrt N \frac{L}{\alpha}\pran{\log{\frac{R_0}{\epsilon}}}^3}\ .$$

    In contrast, the total oracle cost of deterministic restarted methods solving primal-dual problem with global sharpness \cite{applegate2021faster} (which is the optimal deterministic algorithm) is $$\mathcal{O}\pran{N\frac{L_F}{\alpha}\log{\frac{R_0}{\epsilon}}}\ ,$$
    where $L_F$ is the Lipschitz constant for $\mathcal L(x,y)$.
    In the context of finite-sum model, $L_F$ is generally of the same order of $L$ (upto a variance term that is independent of $N$). This shows that RsEGM in general needs $O(\sqrt{N})$ less stochastic oracles to find an $\epsilon$-solution than its deterministic counterpart after ignoring the $\log$ terms. 
\end{rem}

We know from Example \ref{exam:bd} and Example \ref{exam:bilinear} that bilinear problems with bounded region and unconstrained bilinear problem are both globally sharp. A direct application of Theorem \ref{thm:thm-global} gives:
\begin{cor}
    For bilinear problems with bounded region and for unconstrained bilinear problem, the total number of iteration to get $\epsilon$ accuracy with probability at least $1-\delta$ is of order $$\widetilde{\mathcal{O}}\pran{\frac{1}{\sqrt p}\frac{L}{\alpha}\frac{1}{\delta^2}\pran{\log\frac{R_0}{\epsilon}}^3}\ .$$ 
\end{cor}

To prove Theorem \ref{thm:thm-global}, we first introduce two properties of sEGM.
The next lemma bounds the expected distance between the iterates of sEGM and initial point by the distance between initial point and optimal solution. 
\begin{lem}\label{cor:bound-dist}
Consider sEGM (Algorithm \ref{alg:segm}) for solving \eqref{eq:poi}. Suppose Assumption \ref{ass:problem} and Assumption \ref{ass:oracal} holds, then
    \begin{equation}\label{eq:bound-dist}
        \mathbb{E}[\Vert z^{k+1/2}-z^0\Vert] \leq \pran{3+\sqrt{\frac 2 {1-p}}}\Vert z^0-z_* \Vert\ .
    \end{equation}
\end{lem}

\begin{proof}
It follows from Lemma \ref{lem:lem1}
\begin{equation*}
    \mathbb{E}[\Vert z^{k+1/2}-z^{k+1}\Vert^2]\leq 2\phi_0(z_*) \leq 4 \Vert z^0-z_* \Vert^2\ ,
\end{equation*}
thus
\begin{equation*}
        \mathbb{E}[\Vert z^{k+1/2}-z^{k+1}\Vert]\leq \pran{\mathbb{E}[\Vert z^{k+1/2}-z^{k+1}\Vert^2]}^{1/2}\leq 2 \Vert z^0-z_* \Vert\ .
\end{equation*}
We can also derive from the monotonicity of $\mathbb{E}[\phi_{k}(z_*)]$ that
{\footnotesize \begin{equation*}
        \mathbb{E}[\Vert z^{k+1}-z_*\Vert ^2]\leq \frac{1}{1-p} \mathbb{E}[\phi_{k+1}(z_*)]\leq \cdots \leq \frac{1}{1-p}\mathbb{E}[\phi_{0}(z_*)]= \frac{1}{1-p}(2-p)\Vert z^0-z_* \Vert^2 \leq \frac{2}{1-p}\Vert z^0-z_* \Vert^2
\end{equation*}}
thus
\begin{equation*}
        \mathbb{E}[\Vert z^{k+1}-z_*\Vert]\leq \pran{\mathbb{E}[\Vert z^{k+1}-z_*\Vert^2]}^{1/2} \leq \sqrt{\frac{2}{1-p}}\Vert z^0-z_* \Vert\ .
\end{equation*}
By triangle inequality we have
\begin{equation*}
        \mathbb{E}[\Vert z^{k+1/2}-z_*\Vert]\leq \mathbb{E}[\Vert z^{k+1/2}-z^{k+1}\Vert ]+ \mathbb{E}[\Vert z^{k+1}-z_*\Vert] \leq \pran{2+\sqrt{\frac 2 {1-p}}}\Vert z^0-z_* \Vert\ ,
\end{equation*}
\begin{equation*}
        \mathbb{E}[\Vert z^{k+1/2}-z^0\Vert]\leq \mathbb{E}[\Vert z^{k+1/2}-z_*\Vert]+\Vert z_{*}-z^0\Vert \leq \pran{3+\sqrt{\frac 2 {1-p}}}\Vert z^0-z_* \Vert\ .
\end{equation*}
\
\end{proof}

Next, we present the sublinear rate of sEGM on localized duality gap:
\begin{lem}\label{cor:cor-thm1}
Consider sEGM (Algorithm \ref{alg:segm}) for solving \eqref{eq:poi}. Suppose Assumption \ref{ass:problem} and Assumption \ref{ass:oracal} holds, then
    \begin{equation*}
    \mathbb E \left[\max_{z\in \mathcal C\cap\mathcal Z} \mathcal L(\widetilde x^K,y)-\mathcal L(x,\widetilde y^K)\right] \leq \frac 1 {2\tau K} \pran{\frac 72 \max_{z\in \mathcal C}\Vert z^0-z\Vert^2+14\Vert z^0-z_*\Vert^2}\ .
    \end{equation*}
\end{lem}
\begin{proof}
Notice that
\small{\begin{align*}
        \begin{split}
        & \ \mathbb E \left[\max_{z\in \mathcal C\cap\mathcal Z} \mathcal L(\widetilde x^K,y)-\mathcal L(x,\widetilde y^K)\right] \leq \mathbb E \left[\max_{z\in \mathcal C\cap\mathcal Z} \frac 1k \sum_{k=0}^{K-1}\pran{\mathcal L(x^{k+1/2},y)-\mathcal L(x,y^{k+1/2})}\right]\\
        & \ \leq \frac 1k \mathbb E \left[\max_{z\in \mathcal C\cap\mathcal Z} \sum_{k=0}^{K-1} \nabla_x \mathcal L(x^{k+1/2},y^{k+1/2}) (x^{k+1/2}-x)-\nabla_y \mathcal L(x^{k+1/2},y^{k+1/2})(y^{k+1/2}-y)\right]\\
        & \ \leq \frac 1k \mathbb{E}\left[ \max_{z\in \mathcal C}\sum_{k=0}^{K-1} \Theta_{k+1/2}(z)\right] \leq \frac 1 {2\tau K} \pran{\frac 72 \max_{z\in \mathcal C}\Vert z^0-z\Vert^2+14\Vert z^0-z_*\Vert^2}\ .
        \end{split}
\end{align*}}
where the first two inequalities are from the convexity-concavity of $L(x,y)$ and the last inequality is due to Theorem \ref{thm:thm1}.
\end{proof}

{Now we are ready to prove Theorem \ref{thm:thm-global}.}

\begin{proof}[Proof of Theorem \ref{thm:thm-global}]
{Consider the $t$-th outer iteration, and define $R_t=\text{dist}(z^{t,0},\mathcal Z^*)$. To prove the theorem, we'll show that the event $\bigg \{R_T\le \frac{R_0}{2^T}\bigg \}$ is a high probability event and the key is to use the global sharpness on a carefully controlled region.  }

Suppose $K\ge 2000\ge \exp(e^2)$. It follows from Markov inequality that for any inner iteration $k$
\begin{equation*}
    \mathbb P \pran{\Vert z^{t,k}-z^{t,0} \Vert \leq r_0} =1-\mathbb P \pran{\Vert z^{t,k}-z^{t,0} \Vert > r_0} \geq 1-\frac {\mathbb E\left[\Vert z^{t,k}-z^{t,0} \Vert\right]} {r_0}\ .
\end{equation*}
Furthermore, it follows from Lemma \ref{cor:bound-dist} that
\begin{equation}\label{eq:prob-in-region}
    \mathbb P \pran{\Vert z^{t,k}-z^{t,0} \Vert \leq r_0} \geq 1-\frac {\mathbb E[\Vert z^{t,k}-z^{t,0} \Vert]} {r_0} \geq 1-\pran{3+\sqrt{\frac 2 {1-p}}}\frac{\Vert z^{t,0}-z_* \Vert}{r_0}\ .
\end{equation}
Thus, if $r_0$ is chosen to be sufficiently large, the probability $\mathbb P \pran{\Vert z^{t,k}-z^{t,0} \Vert \leq r_0}$ can be arbitrarily close to $1$ . Given the high probability event $\left\{ \Vert z^{t,k}-z^{t,0} \Vert \leq r_0 \right\}$, we can upper bound the normalized duality gap as 
\begin{equation}\label{eq:eq-dgap-upperbd}
     \max_{z\in W_{r}(z^{t,k})}{\left\{\mathcal L(x^{t,k},y)-\mathcal L(x,y^{t,k})\right\}} \leq  \max_{z\in W_{r+r_0}(z^{t,0})}{\left\{\mathcal L(x^{t,k},y)-\mathcal L(x,y^{t,k})\right\}}\ .
\end{equation}
By $\alpha$-sharpness and the fact that $W_{r+r_0}(z^{t,0})$ is deterministic given $z^{t,0}$, we have for any constant $C> 0$ and function $g(k)> 0$ that
\begin{align*}
    \begin{split}
        \mathbb P \pran{\alpha R_{t+1}\leq \frac{CR_t}{g(k)}} \geq & \ \mathbb P \pran{\frac 1r \max_{z\in W_r(z^{t+1,0})}{\mathcal L(x^{t+1,0},y)-\mathcal L(x,y^{t+1,0})}\leq \frac{CR_t}{g(k)}}\\
        = & \ \mathbb P \pran{\max_{z\in W_r(z^{t,k})}{\mathcal L(x^{t,k},y)-\mathcal L(x,y^{t,k})}\leq \frac{CrR_t}{g(k)}} \ .
    \end{split}
\end{align*}

Thus, it holds for any $r_0\ge 0$ that
{\footnotesize{\begin{align}\label{eq:main-bounded}
    \begin{split}
    \mathbb P \pran{\alpha R_{t+1}\leq \frac{CR_t}{g(k)}} \geq & \ \mathbb P \pran{\max_{z\in W_r(z^{t,k})}{\mathcal L(x^{t,k},y)-\mathcal L(x,y^{t,k})}\leq \frac{CrR_t}{g(k)} \;, \; \Vert z^{t,k}-z^{t,0} \Vert \leq r_0}\\
    \geq & \ \mathbb P \pran{\max_{z\in W_{r+r_0}(z^{t,0})}{\mathcal L(x^{t,k},y)-\mathcal L(x,y^{t,k})}\leq \frac{CrR_t}{g(k)} \;, \; \Vert z^{t,k}-z^{t,0} \Vert \leq r_0}\\
    = & \ \mathbb P \pran{\Vert z^{t,k}-z^{t,0} \Vert \leq r_0} \\
    & \ - \mathbb P \pran{\max_{z\in W_{r+r_0}(z^{t,0})}{\mathcal L(x^{t,k},y)-\mathcal L(x,y^{t,k})} > \frac{CrR_t}{g(k)} \;, \;  \Vert z^{t,k}-z^{t,0} \Vert \leq r_0}\\
    \geq & \ \mathbb P \pran{\Vert z^{t,k}-z^{t,0} \Vert \leq r_0} - \mathbb P \pran{\max_{z\in W_{r+r_0}(z^{t,0})}{\mathcal L(x^{t,k},y)-\mathcal L(x,y^{t,k})} > \frac{CrR_t}{g(k)}} \\
    \end{split}
\end{align}}}
where the second inequality utilizes Equation \eqref{eq:eq-dgap-upperbd}.

The next step is to upper bound the last term in the right hand side of \eqref{eq:main-bounded}. Denote $C_1=3+\sqrt{\frac 2 {1-p}}$ and choose $r_0=r=2C_1R_t\frac{T}{\delta}$, then we have $z_*\in W_{r+r_0}(z^{t,0})$ because $r+r_0\geq R_t$, thus
{\footnotesize \begin{align}\label{eq:eq-dgap-markov}
    \begin{split}
        \mathbb P & \pran{\max_{z\in W_{r+r_0}(z^{t,0})}  {\mathcal L(x^{t,k},y)-\mathcal L(x,y^{t,k})} > \frac{CrR_t}{g(k)} } \leq \frac{g(k)}{CrR_t} \mathbb E \left [\max_{z\in W_{r+r_0}(z^{t,0})}{\mathcal L(x^{t,k},y)-\mathcal L(x,y^{t,k})} \right ]\\
        & \ \leq \frac{g(k)}{CrR_t} \frac 1 {2\tau k} \pran{\frac 72+14} \max_{z\in W_{r+r_0}(z^{t,0})}\Vert z^{t,0}-z\Vert^2  \ \leq \frac{35}{4\tau C}\frac{g(k)}{k}\frac{(r+r_0)^2}{rR_t} \ ,
    \end{split}
\end{align}}
where the first inequality follows from Markov inequality, the second one utilizes Lemma \ref{cor:cor-thm1} and $z_*\in W_{r+r_0}(z^{t,0})$.
Combine Equation \eqref{eq:main-bounded} with \eqref{eq:eq-dgap-markov}, we obtain
\begin{align*}
    \begin{split}
        \mathbb P \pran{\alpha R_{t+1}\leq \frac{CR_t}{g(k)}} & \ \geq \mathbb P \pran{\Vert z^{t,k}-z^{t,0} \Vert \leq r_0} - \frac{35}{4\tau C}\frac{g(k)}{k}\frac{(r+r_0)^2}{rR_t}\\
        & \ \geq 1-\pran{3+\sqrt{\frac 2 {1-p}}}\frac{R_t}{r_0} - \frac{35}{4\tau C}\frac{g(k)}{k}\frac{(r+r_0)^2}{rR_t}\ ,
    \end{split}
\end{align*}
where the last inequality uses Equation \eqref{eq:prob-in-region}.
Thus, we obtain by substituting the value of $r_0$, $r$ and $\tau$ that
\begin{align*}
    \begin{split}
    \mathbb P \pran{\alpha R_{t+1}\leq \frac{CR_t}{g(k)}} \geq & 1-\frac{\delta}{2T} - \frac{35}{4\tau C}\frac{g(k)}{k}\frac{8C_1T}{\delta} 
    =  1-\frac{\delta}{2T}-\frac{L}{C}\frac{g(k)}{k}\frac{T}{\delta}\frac{140C_1}{\sqrt{p}}\ .
    \end{split}
\end{align*}
Let $C_2=140C_1$, $C=\alpha$ and function $g(k)=\log{\log{k}}$, then we have 
\begin{equation}\label{eq:ieq-3}
    \begin{split}
    \mathbb P \pran{R_{t+1}\leq \frac{R_t}{\log{\log{k}}}}
    \geq 1-\frac{\delta}{2T}-\frac{C_2}{\sqrt p}\frac{L}{\alpha}\frac{T}{\delta}\frac{\log{\log{k}}}{k} \ .
    \end{split}
\end{equation}

When $k\geq 2000 > \exp{(e^2)}$, we have $\left \{R_{t}\leq \frac{R_{t-1}}{\log{\log{k}}} \right \} \subseteq \left \{R_{t}\leq \frac {R_{t-1}}{2} \right \}$. Let $\mathcal E_t=\left\{R_{t}\leq \frac{R_{t-1}}{2}\right\}$ and $\mathbb P\pran{\cdot|\bar{\mathcal{E}}_t}$ be the probability conditioned on $\left\{\mathcal E_1,...,\mathcal E_t\right\}$, then 
\begin{align}\label{eq:prob1}
    \begin{split}
    \mathbb P\pran{R_T\leq \frac {R_0}{2^T} } \geq & \ \mathbb P\pran{\cap_{t=1}^{T}\mathcal E_t} = \mathbb P(\mathcal E_1)\prod_{t=2}^{T}\mathbb P(\mathcal E_t | \bar{\mathcal{E}}_{t-1})\\
    \geq & \ \bigg( 1-\frac{\delta}{2T}-\frac{C_2}{\sqrt p}\frac{L}{\alpha}\frac{T}{\delta}\frac{\log{\log{k}}}{k}\bigg )^T
    \geq  \ 1-\frac{\delta}{2}-\frac{C_2}{\sqrt p}\frac{L}{\alpha}\frac{T^2}{\delta}\frac{\log{\log{k}}}{k}\ ,
    \end{split}
\end{align}
where the second inequality is from Equation \eqref{eq:ieq-3} and the third one uses Bernoulli inequality.
Thus, if $k\geq \frac{2C_2}{\sqrt p}\frac{L}{\alpha}\frac{T^2}{\delta^2}\log{\pran{\frac{2C_2}{\sqrt p}\frac{L}{\alpha}\frac{T^2}{\delta^2}}}$ we obtain
\begin{align}\label{eq:prob2}
    \begin{split}
        \frac{C_2}{\sqrt p}\frac{L}{\alpha}\frac{T^2}{\delta}\frac{\log{\log{k}}}{k} \leq & \frac{\delta}{2} \cdot \frac{\log{\log{\pran{\frac{2C_2}{\sqrt p}\frac{L}{\alpha}\frac{T^2}{\delta^2}\log{\pran{\frac{2C_2}{\sqrt p}\frac{L}{\alpha}\frac{T^2}{\delta^2}}}}}}}{\log{\pran{\frac{2C_2}{\sqrt p}\frac{L}{\alpha}\frac{T^2}{\delta^2}}}}\\
        \leq & \frac{\delta}{2} \cdot \frac{\log{\pran{\log{\pran{\frac{2C_2}{\sqrt p}\frac{L}{\alpha}\frac{T^2}{\delta^2}}+\log{\log{\pran{\frac{2C_2}{\sqrt p}\frac{L}{\alpha}\frac{T^2}{\delta^2}}}}}}}}{\log{\pran{\frac{2C_2}{\sqrt p}\frac{L}{\alpha}\frac{T^2}{\delta^2}}}} 
        \leq  \frac{\delta}{2}\ ,
    \end{split}
\end{align}
where the last inequality is from $\frac{\log{(x+\log{x})}}{x}\leq 1$.
Combine \eqref{eq:prob1} and \eqref{eq:prob2} we arrive at
\begin{equation*}
    \mathbb P\pran{R_T\leq \frac {R_0}{2^T}} \geq 1-\frac{\delta}{2}-\frac{C_2}{\sqrt p}\frac{L}{\alpha}\frac{T^2}{\delta}\frac{\log{\log{k}}}{k} \geq 1-\delta\ .
\end{equation*}

Finally, by choosing $T=\left\lceil \log_{2}\frac{R_0}{\epsilon}\right\rceil$, $K\geq \max\left\{\frac{2C_2}{\sqrt p}\frac{L}{\alpha}\frac{T^2}{\delta^2}\log{\pran{\frac{2C_2}{\sqrt p}\frac{L}{\alpha}\frac{T^2}{\delta^2}}},2000 \right\}$, the total number of iteration to get $\epsilon$ accuracy with probability at least $1-\delta$ is $TK$ which is of order 
\begin{equation*}
\widetilde{\mathcal{O}}\pran{\frac{1}{\sqrt p}\frac{L}{\alpha}\frac{1}{\delta^2}\pran{\log\frac{R_0}{\epsilon}}^3}\ .
\end{equation*}
\ 
\end{proof}

\subsection{Convergence Guarantees for Standard-Form LP}

In this subsection, we establish the high probability convergence rate of RsEGM for LP \eqref{eq:lp}.  Unfortunately, LP is not globally sharp (see Appendix \ref{sec:lp_not_global_sharp}), thus the result in Section \ref{sec:global-sharp} cannot be directly applied to LP. Instead, \cite{applegate2021faster} shows that LP is sharp on any bounded region. Unlike deterministic methods as studied in \cite{applegate2021faster}, where the iterates of a primal dual algorithm always stay in a bounded region, the iterates of stochastic methods 
may escape from any bounded region. Theorem \ref{thm:thm-LP} overcomes this difficulty by presenting a high probability argument.


\begin{thm}\label{thm:thm-LP}
    Consider RsEGM (Algorithm \ref{alg:rsegm}) for solving LP \eqref{eq:lp}. Suppose Assumption \ref{ass:problem} and Assumption \ref{ass:oracal} hold. Let $L$ be the Lipschitz parameter of the stochastic oracle and $H$ be the Hoffman constant of the KKT system of LP (see Example \ref{thm:lem-LP-sharp}). Denote the initial distance to optimal set as $R_0=\text{dist}(z^{0,0},\mathcal Z^*)$, and choose step-size $\tau=\frac{\sqrt{p}}{2L}$. For any $\delta\in(0,1)$  and $\epsilon>0$, if the outer loop count $T$ and the inner loop count $K$ of RsEGM satisfy $$T\ge \max{\left\{ \log_{2}\frac{R_0}{\epsilon}, 0 \right\}},\; K\geq \widetilde{\mathcal{O}}\pran{\frac{1}{\sqrt p}\frac{L}{1/\left[H(1+\Vert z^{0,0} \Vert +R_0)\right]}\frac{1}{\delta^3}T^4},$$ then it holds with probability at least $1-\delta$ that $$\text{dist} (z^{T,0},\mathcal Z^*)\leq \epsilon.$$ Furthermore, the total number of iteration to get $\epsilon$ accuracy with probability at least $1-\delta$ is of order $$\widetilde{\mathcal{O}}\pran{\frac{1}{\sqrt p}\frac{L}{1/[H(1+\Vert z^{0,0} \Vert +R_0)]}\frac{1}{\delta^3}\pran{\log\frac{R_0}{\epsilon}}^5}\ . $$
\end{thm}

\begin{rem}
      Similar to Remark \ref{rem:median}, we can obtain the complexity of the median trajectory of the algorithm by choosing $\delta=\frac{1}{2}$.
\end{rem}

\begin{rem}
    As shown in \cite[Lemma 5]{applegate2021faster}, LP is sharp with $\alpha=\frac{1}{H\sqrt{1+4R^2}}$ in the ball $\{z|\|z\|_2\le R\}$ with any radius $R>0$. Notice that for deterministic primal-dual methods, the iterates usually stay in a ball centered at $0$ with radius $R=\|z^{0,0}\|+R_0$ (see Section 2 in \cite{applegate2021faster}), thus the effective sharpness constant on deterministic algorithms for LP is given by $\alpha=\frac{1}{H(1+\|z^{0,0}\|+R_0)}$. Theorem \ref{thm:thm-LP} shows that the complexity of RsEGM is $\widetilde{\mathcal{O}}\pran{\frac{1}{\sqrt p}\frac{L}{1/\alpha}\frac{1}{\delta^3}\pran{\log\frac{R_0}{\epsilon}}^5}\ . $ 
\end{rem}



{In the following we prove Theorem \ref{thm:thm-LP}. The key is to show that the iterates stay in a bounded region with a high probability, thus we can utilize the sharpness condition of LP.}

\begin{proof}[Proof of Theorem \ref{thm:thm-LP}]

Suppose $K\ge 2000\ge \exp(e^2)$. Denote $R_t=\text{dist}(z^{t,0},\mathcal Z^*)$ and consider event $\mathcal E_t=\left\{R_{t}\leq \frac{1}{\log{\log{K}}}R_{t-1}\; \text{and}\; \Vert z^{t,0}-z^{t-1,0} \Vert \leq R_{t-1}B\right\}$, where we determine the value of $B$ later.  By definition:
\begin{equation*}
    \cap_{l=1}^{t}\mathcal E_l \subseteq \{R_{t}<R_{t-1}<...<R_0 \; \text{and} \; \Vert z^{t,0} \Vert \leq tR_0B +\Vert z^{0,0} \Vert \}
\end{equation*}

Let $\bar{\mathcal{E}}_t=\{\cap_{l=1}^{t}\mathcal E_l\}$, then $\mathbb P(\cdot|\bar{\mathcal{E}}_t)$ is the probability conditioned on $\{\mathcal E_1,...,\mathcal E_t\}$. For convenience, we denote, $\mathbb P(\cdot|\bar{\mathcal{E}}_0)=\mathbb P(\cdot)$.
It follows from Lemma \ref{thm:lem-LP-sharp} and the sharpness of LP that
{\small \begin{align}\label{eq:eq-LP-sharp}
    \begin{split}
        & \ \mathbb P \pran{\mathcal E_{t+1}|\bar{\mathcal{E}}_t} \\  = & \ \mathbb P \pran{\frac{1}{H\sqrt{1+4\Vert z^{t+1,0} \Vert ^2}} R_{t+1}\leq \frac{1}{H\sqrt{1+4\Vert z^{t+1,0} \Vert ^2}} \frac{R_{t}}{\log{\log{K}}} ,\; \Vert z^{t+1,0}-z^{t,0} \Vert \leq R_{t}B \big| \bar{\mathcal{E}}_t}\\
        \geq & \ \mathbb P \pran{\max_{\hat z\in W_r(z^{t+1,0})}\frac {\mathcal L(x^{t+1,0},\hat y)-\mathcal L(\hat x,y^{t+1,0})}{r} \leq \frac{R_{t}}{ H \log{\log{K}}  \sqrt{1+4\Vert z^{t+1,0} \Vert ^2}}, \; \Vert z^{t+1,0}-z^{t,0} \Vert \leq R_{t}B \big| \bar{\mathcal{E}}_t}\ . \\
    \end{split}
\end{align}}
Note that given the event $\{\Vert z^{t+1,0}-z^{t,0} \Vert \leq R_tB\}$, we have for any $r>0$ that
\begin{equation}\label{eq:eq-ball-relax-1}
    \max_{\hat z\in W_r(z^{t+1,0})}\frac {\mathcal L(x^{t+1,0},\hat y)-\mathcal L(\hat x,y^{t+1,0})}{r} \leq \max_{\hat z\in W_{r+R_tB}(z^{t,0})}\frac {\mathcal L(x^{t+1,0},\hat y)-\mathcal L(\hat x,y^{t+1,0})}{r} \ .
\end{equation}
Furthermore, we have conditioned on $\bar{\mathcal{E}}_t$ that $\Vert z^{t,0} \Vert \leq tR_0B+\Vert z^{0,0} \Vert$, thus
\begin{equation}\label{eq:eq-ball-relax-2}
    \frac{R_t}{H \log{\log{K}} \sqrt{1+4(\Vert z^{t,0} \Vert+R_tB) ^2}} \geq \frac{R_t}{H\log{\log{K}}\sqrt{1+4(\Vert z^{0,0} \Vert+t R_0 B+R_0B) ^2}} \ .
\end{equation}

Combining Equation \eqref{eq:eq-LP-sharp}, \eqref{eq:eq-ball-relax-1} and \eqref{eq:eq-ball-relax-2}, we arrive at
{\small
\begin{align}\label{eq:eq-premarkov-LP}
    \begin{split}
        & \mathbb P \pran{\mathcal E_{t+1}|\bar{\mathcal{E}}_t} \\  
        \geq & \mathbb P \pran{\max_{\hat z\in W_{r+R_{t}B}(z^{t,0})}\frac {\mathcal L(x^{t,K},\hat y)-\mathcal L(\hat x,y^{t,K})}{r} \leq \frac{R_t/\log{\log{K}}}{H\sqrt{1+4(\Vert z^{0,0} \Vert+(t+1) R_0 B) ^2}}, \Vert z^{t,K}-z^{t,0} \Vert \leq R_{t}B \big| \bar{\mathcal{E}}_t}\\
        = & \mathbb P \big(\Vert z^{t,K}-z^{t,0} \Vert \leq R_{t}B \big| \bar{\mathcal{E}}_t \big)\\
        - & \mathbb P \pran{\max_{\hat z\in W_{r+R_{t}B}(z^{t,0})}\frac {\mathcal L(x^{t,K},\hat y)-\mathcal L(\hat x,y^{t,K})}{r} > \frac{R_t/\log{\log{K}}}{H\sqrt{1+4(\Vert z^{0,0} \Vert+(t+1) R_0 B) ^2}}, \Vert z^{t,K}-z^{t,0} \Vert \leq R_{t}B \big| \bar{\mathcal{E}}_t}\\
        \geq & \mathbb P \pran{\Vert z^{t,K}-z^{t,0} \Vert \leq R_{t}B}-\mathbb P \pran{\max_{\hat z\in W_{r+R_{t}B}(z^{t,0})}\frac {\mathcal L(x^{t,K},\hat y)-\mathcal L(\hat x,y^{t,K})}{r} > \frac{R_t/\log{\log{K}}}{H\sqrt{1+4(\Vert z^{0,0} \Vert+(t+1) R_0 B) ^2}}\big| \bar{\mathcal{E}}_t}\ .
    \end{split}
\end{align}}

Next, we set $C_1=3+\sqrt{\frac 2 {1-p}}$, $B=\frac{2TC_1}{\delta}$ and $r=R_t B$, then we know there exists an optimal solution $z_*$ such that $z_*\in W_{r+R_tB}(z^{t,0})$, because $r+R_tB=2R_tB>R_t$. Then we can upper bound the second term in the right hand side of \eqref{eq:eq-premarkov-LP}  by Markov inequality:
{\small \begin{align}\label{eq:eq-markov-LP}
    \begin{split}
        & \ \mathbb P \pran{\max_{\hat z\in W_{r+R_{t}B}(z^{t,0})}\frac {\mathcal L(x^{t,K},\hat y)-\mathcal L(\hat x,y^{t,K})}{r} > \frac{R_t/\log{\log{K}}}{H\sqrt{1+4(\Vert z^{0,0} \Vert+(t+1) R_0 B) ^2}}\big| \bar{\mathcal{E}}_t}\\
        & \ \leq \frac{H \log{\log{K}}\sqrt{1+4(\Vert z^{0,0} \Vert+(t+1) R_0 B) ^2}}{rR_t}\mathbb E\left[\max_{\hat z\in W_{r+R_{t}B}(z^{t,0})}\mathcal L(x^{t,K},\hat y)-\mathcal L(\hat x,y^{t,K})\right] \\
        & \ \leq \frac{H\log{\log{K}}\sqrt{1+4(\Vert z^{0,0} \Vert+(t+1) R_0 B) ^2}}{rR_t} \frac {35} {4\tau K} max_{z\in W_{r+R_{t}B}(z^{t,0})}\Vert z^{t,0}-z\Vert^2\\
        & \ \leq \frac{H\log{\log{K}}\sqrt{1+4(\Vert z^{0,0} \Vert+(t+1) R_0 B) ^2}}{rR_t} \frac {35} {4\tau K} (r+R_{t}B)^2\\
        & \ \leq \frac{H\log{\log{K}}(1+2(\Vert z^{0,0} \Vert+(t+1) R_0 B))}{rR_t} \frac {35} {4\tau K} (r+R_{t}B)^2 \ ,
    \end{split}
\end{align}}
where the second inequality follows from Lemma \ref{cor:cor-thm1}, by noticing $z_*\in W_{r+R_tB}(z^{t,0})$, and the last inequality uses $\sqrt{1+x}\leq 1+\sqrt{x}$. 
Combining Equation \eqref{eq:prob-in-region}, \eqref{eq:eq-premarkov-LP} and \eqref{eq:eq-markov-LP}, we obtain
{\small \begin{equation*}
    \mathbb P (\mathcal E_{t+1}|\bar{\mathcal{E}}_t) \geq 1-\pran{3+\sqrt{\frac 2 {1-p}}}\frac{1}{B}-\frac {35H} {4\tau}\frac {\log{\log{K}}}{K}\frac{(r+R_{t}B)^2}{rR_t}\pran{1+2(\Vert z^{0,0} \Vert+(t+1) R_0 B)} 
\end{equation*}}

Furthermore, by recalling $r=R_tB$, we obtain
\begin{align*}
    \begin{split}
    &\mathbb P\pran{R_T\leq \frac{R_0}{2^T}}\\
    \geq & \ \mathbb P(\cap_{t=1}^{T}\mathcal E_t) = \mathbb P(\mathcal E_1)\prod_{t=2}^{T}\mathbb P(\mathcal E_t | \bar{\mathcal{E}}_{t-1})\\
    \geq & \ 1-\pran{3+\sqrt{\frac 2 {1-p}}}\frac{T}{B} \ -\frac {35H} {4\tau} \frac{\log{\log{K}}}{K}\frac{(B+B)^2}{B}\pran{T+2T\Vert z^{0,0} \Vert+T(T+1)R_0B} \\
    \end{split}
\end{align*}
where the last inequality comes from $\prod_i (1-x_i)\geq 1-\sum_i x_i$ for $x_i>0$.
Thus we have
{\small \begin{equation*}
    \mathbb P\pran{R_T\leq \frac{R_0}{2^T}} \geq 1-\pran{3+\sqrt{\frac 2 {1-p}}}\frac{T}{B} -\frac {35HB} {\tau} \frac{\log{\log{K}}}{K}\pran{T+2T\Vert z^{0,0} \Vert+T(T+1)R_0B} 
\end{equation*}}

Recall $\tau=\frac{\sqrt{p}}{2L}$ and $B=\frac{2TC_1}{\delta}$, thus
\begin{align*}
\begin{split}
    \mathbb P\pran{R_T\leq \frac{R_0}{2^T}} \geq & \ 1-\frac{\delta}{2} -\frac {140HL}{\sqrt{p}} \frac{\log{\log{K}}}{K}\frac{T}{\delta}\pran{T+2T\Vert z^{0,0} \Vert+T(T+1)R_0\frac{2TC_1}{\delta}}\\
    \geq & \ 1-\frac{\delta}{2} -\frac {140HL}{\sqrt{p}}\frac{\log{\log{K}}}{K}\frac{T}{\delta}\frac{T^3}{\delta}(1+2\Vert z^{0,0} \Vert+4C_1R_0)\\
    \geq & \ 1-\frac{\delta}{2} - HL\frac{\log{\log{K}}}{K}\frac{T^4}{\delta^2}\frac{140}{\sqrt{p}}(1+2\Vert z^{0,0} \Vert+4C_1R_0) \ .
\end{split}
\end{align*}

Denote $C_2=C_2(z^{0,0})=140(1+2\Vert z^{0,0} \Vert+4C_1R_0)$. Then if $$K\geq \frac{2C_2}{\sqrt p}HL\frac{T^4}{\delta^3}\log{\pran{\frac{2C_2}{\sqrt p}HL\frac{T^4}{\delta^3}}},$$ we have $$HL\frac{\log{\log{K}}}{K}\frac{T^4}{\delta^2}\frac{140}{\sqrt{p}}(1+2\Vert z^{0,0} \Vert+4C_1R_0) \leq \frac{\delta}{2}\ , $$
for the same reason as Equation \eqref{eq:prob2}.
Thus, we obtain
\begin{equation*}
    \mathbb P\pran{R_T\leq \frac{R_0}{2^T}} \geq 1-\frac{\delta}{2} - HL\frac{\log{\log{K}}}{K}\frac{T^4}{\delta^2}\frac{140}{\sqrt{p}}(1+2\Vert z^{0,0} \Vert+4C_1R_0) \geq 1-\delta \ .
\end{equation*}

In summary, choosing $T=\left\lceil \log_{2}\frac{R_0}{\epsilon}\right\rceil$, $K\geq \max\left\{\frac{2C_2}{\sqrt p}HL\frac{T^4}{\delta^3}\log{\pran{\frac{2C_2}{\sqrt p}HL\frac{T^4}{\delta^3}}}, 2000\right\}$, the total number of iteration to get $\epsilon$ accuracy with probability at least $1-\delta$ is $TK$ which is of order \begin{equation*}
\widetilde{\mathcal{O}}\pran{\frac{1}{\sqrt p}\frac{L}{1/[H(1+\Vert z^{0,0} \Vert +R_0)]}\frac{1}{\delta^3}\pran{\log\frac{R_0}{\epsilon}}^5}\ .
\end{equation*}
\ 
\end{proof}

\section{Lower Bound}\label{sec:lower-bound}
In this section, we establish the complexity lower bound of stochastic first-order algorithms (with variance reduction) for solving sharp primal-dual problems. Together with the results shown in Section \ref{sec:convergence}
, we demonstrate that RsEGM has nearly optimal convergence rate (upto log terms) among a large class of stochastic primal-dual methods.

Since we study the lower bound in this section, we can assume $\mathcal{L}$ is differentiable everywhere, thus $\nabla \mathcal{L}$ is well defined. Recall that in deterministic first-order primal-dual methods, we utilize the gradient information to update the iterates. Thus, the iterates of first-order methods always stay in the subspace spanned by the gradient, and we call this behavior first-order span-respecting (see for example \cite{applegate2021faster}
) defined as below:
\begin{mydef}\label{def:def-span}
    A deterministic primal dual algorithm is called first-order  span-respecting if its iterates satisfy
    \begin{align*}
        \begin{split}
            & \ x^t\in x^0+\text{Span}\left\{\nabla_x\mathcal L(x^i,y^j):\;\forall i,j\in\{1,...,t-1\}\right\}\\
            & \ y^t \in y^0+\text{Span}\left\{\nabla_y\mathcal L(x^i,y^j):\;\forall i\in\{1,...,t\},\forall j\in\{1,...,t-1\}\right\}
        \end{split}
    \end{align*}
\end{mydef}

Many classic primal dual first-order methods, such as GDA, AGDA and EGM, are first-order span-respecting when applying unconstrained problems.

In this paper, we study stochastic algorithms with variance reduction. As a result, we consider a general class of stochastic algorithms where the span includes both the deterministic gradient (for variance reduction sake) and the stochastic oracle $F_\xi$.  More formally, we introduce stochastic span-respecting algorithms, defined as below:


\begin{mydef}\label{def:def-span-random}
    A randomized primal dual algorithm is called  first-order stochastic span-respecting with respect to a stochastic oracle $\xi$ if the iterates in $\mathbb R^m$ satisfy:
    \begin{align*}
        \begin{split}
            & \ x^t\in x^0+\text{Span}\left\{\nabla_x\mathcal L(x^i,y^j),\nabla_{x}^{\xi^{i,j}}\mathcal L(x^i,y^j):\;\forall i,j\in\{1,...,t-1\}\right\}\\
            & \ y^t \in y^0+\text{Span}\left\{\nabla_y\mathcal L(x^i,y^j),\nabla_{y}^{\xi^{i,j}}\mathcal L(x^i,y^j):\;\forall i\in\{1,...,t-1\},\forall j\in\{1,...,t\}\right\}
        \end{split}
    \end{align*}
     $\xi^{i,j}$ are random variables.
\end{mydef}

\begin{rem}
    The above definition consists of both deterministic gradient together with the stochastic gradient estimator. The reason is to include the variance reduction method into the algorithm class. Obviously the above algorithm class is larger than the one only with deterministic or stochastic gradient. As a result, the lower bound is generally higher than a pure stochastic algorithm.
\end{rem}



\begin{rem}
    If $\mathcal L$ is bilinear, then with appropriate indexing of iterates, sEGM and RsEGM are first-order stochastic span-respecting.
\end{rem}

\begin{thm}\label{thm:lower}
    Consider any iteration $t\in \mathbb N$ and parameter value $L>\sqrt 2\alpha>0$. There exists an $\alpha$-sharp primal-dual problem  that satisfies Assumption \ref{ass:problem} 
     and an $L$-Lipschitz (in expectation) stochastic oracle $F_{\xi}$ that satisfies Assumption \ref{ass:oracal},
      such that the iterates $z^t$ of any stochastic span-respecting algorithm  satisfies that
    \begin{equation*}\label{eq:lower-bound}
        \text{dist}(z^t,\mathcal Z^*)\geq \frac{1}{\sqrt 2}\pran{1-\frac{\alpha\sqrt{2}}{L}}^{t} \text{dist}(z^0,\mathcal Z^*) \ .
    \end{equation*}
\end{thm}

\begin{rem}
    Theorem \ref{thm:lower} implies the following lower complexity bound for a stochastic first-order method to achieve an $\epsilon$-accuracy solution: $$\Omega\pran{\frac{L}{\alpha}\log{\frac{R_0}{\epsilon}}} \ .$$
\end{rem}

\begin{proof}[Proof of Theorem \ref{thm:lower}]
Set $m=2t$. We consider a primal-dual problem of the form 
\begin{equation}\label{eq:bilinear-problem}
        \min_{x\in\mathbb R^{2m}}\max_{y\in \mathbb R^{2m}} \mathcal L(x,y)=y^TAx-b^Ty \ ,
\end{equation} 
The matrix $A\in \mathbb R^{2m\times 2m}$ and vector $b\in \mathbb R^{2m}$ have the structure
\begin{align*}
    A=\begin{pmatrix}
        A_0 & \ \\ \ & A_0
    \end{pmatrix},\;
    b=\begin{pmatrix}
        b_0 \\
        b_0
    \end{pmatrix} \ ,
\end{align*}
where $A_0\in \mathbb R^{m\times m}$ and $b_0 \in \mathbb R^m$ are determined later.
Rewriting the problem using $x=(x_1,x_2)$ and $y=(y_1,y_2)$, we obtain:
\begin{equation*}
    \min_{x\in\mathbb R^{2m}}\max_{y\in \mathbb R^{2m}} \mathcal L(x,y)=\mathcal L_1(x,y) +\mathcal L_2(x,y) \ ,
\end{equation*} 
where $\mathcal L_1(x,y)=y_1^TA_0x_1-b_0^Ty_1$, $\mathcal L_2(x,y)=y_2^TA_0x_2-b_0^Ty_2$ .
    
For any given Lipschitz parameter $L$ and sharpness parameter $\alpha$ such that $L>\sqrt{2}\alpha$, we denote $Q=\frac{L^2}{2\alpha^2}$, $\kappa=\frac{\sqrt{Q}+3}{\sqrt{Q}+1}$, and 
\begin{align*}
G_0=
\begin{pmatrix}
    2 & -1 & \ & \ & \ & \ \\
    -1 & 2 & -1 & \ & \ & \ \\
    \ & -1 & 2 & -1 & \ & \  \\
    \ & \ & \ddots & \ddots & \ddots & \ \\
    \ & \ & \ & -1 & 2 & -1\\
    \ & \ & \ & \ & -1 & \kappa
\end{pmatrix}\ .
\end{align*}
Note that it holds for any $\kappa\leq 3$ that $0\preceq G_0 \preceq 4I.$
Denote 
\begin{equation*}
    H_0=\frac{L^2/2-\alpha^2}{4}G_0+\alpha^2 I,\; h_0=\frac{L^2/2-\alpha^2}{4}e_1\ ,
\end{equation*}
where $e_1$ is the first standard unit vector. One can check that the solution of the linear system $H_0u=h_0$ is unique and given by $$u^*\in \mathbb R^m \;\text{such that}\; u_i^*=q^i,\; i=1,2,...,m \ ,$$ with $q=\frac{\sqrt{Q}-1}{\sqrt{Q}+1}$. Furthermore, it holds that $$\alpha^2 I \preceq H_0  \preceq \frac{L^2}{2}I \ .$$
Choose $A_0=H_0^{1/2}$, $b_0=H_0^{-1/2}h_0$ and consider \eqref{eq:bilinear-problem}. The optimal solution to \eqref{eq:bilinear-problem} is given by 
    $Ax=b,\; A^Ty=0 \ ,$
which is equivalent to 
\begin{equation*}
    A_0x_1=b_0,\; A_0^Ty_1=0,\; A_0x_2=b_0,\; A_0^Ty_2=0 \ ,
\end{equation*}
and thus
\begin{equation*}
    H_0x_1=h_0,\; y_1=0,\; H_0x_2=h_0,\; y_2=0 \ .
\end{equation*}
The optimal solution to \eqref{eq:bilinear-problem} is given by $z^*=(x_1^*,x_2^*,y_1^*,y_2^*)=(u^*,u^*,0,0)$. 
Moreover, the primal-dual problem \eqref{eq:bilinear-problem} is $\alpha$-sharp, 
because
\begin{equation*}
    \alpha I_m \preceq A_0= H_0^{1/2} \ ,
    \;\; \alpha I_{2m} \preceq A= \begin{pmatrix}
        A_0 & \ \\ \ & A_0
    \end{pmatrix} \ .
\end{equation*}

Without loss of generality we assume $x^0=y^0=0$. 
Consider the stochastic oracle
\begin{align}\label{eq:eq-lower-bound-oracle}
    \begin{split}
    F_\xi (z)& \ =2(\nabla_x \mathcal L_\xi(x,y),-\nabla_y \mathcal L_\xi(x,y)),\\  
             & \ =  \begin{cases} 
                        2(A_0^Ty_1,0,A_0x_1-b,0),\; \text{with probability } \frac 12\\
                        2(0,A_0^Ty_2,0,A_0x_2-b),\; \text{with probability } \frac 12
                    \end{cases}.
    \end{split}
\end{align}
Then we know that $F_{\xi}(z)$ is unbliased and $\mathbb E [\Vert F_{\xi}(u)-F_{\xi}(v) \Vert ^2]\leq L^2\Vert u-v \Vert ^2$, thus the stochastic gradient oracle satisfies Assumption \ref{ass:oracal}.

By definition of span-respecting for randomized algorithms with the above oracle:
\begin{align*}
    \begin{split}
        & \ x_1^t\in x_1^0+ \text{Span}\left\{A_0^Ty_1^0,...,A_0^Ty_1^{t-1}\right\} \\
        & \ x_2^t\in x_2^0+ \text{Span}\left\{A_0^Ty_2^0,...,A_0^Ty_2^{t-1}\right\} \\
        & \ y_1^t\in y_1^0+ \text{Span}\left\{A_0x_1^0-b_0,...,A_0x_1^{t-1}-b_0\right\}\\
        & \ y_2^t\in y_2^0+ \text{Span}\left\{A_0x_2^0-b_0,...,A_0x_2^{t-1}-b_0\right\}
    \end{split}
\end{align*}
thus it holds that
{\scriptsize \begin{align*}
    \begin{split}
        x_1^t & \ \in \text{Span}\left\{A_0^T(A_0x_1^0-b),...,A_0^T(A_0x_1^{t-1}-b)\right\}= \text{Span}\left\{H_0x_1^0-h_0,...,H_0x_1^{t-1}-h_0\right\} \subseteq \{e_1,...,e_t \}\\
        x_2^t & \ \in \text{Span}\left\{A_0^T(A_0x_2^0-b),...,A_0^T(A_0x_2^{t-1}-b)\right\}= \text{Span}\left\{H_0x_2^0-h_0,...,H_0x_2^{t-1}-h_0\right\} \subseteq \{e_1,...,e_t \}
    \end{split}
\end{align*}}

Therefore, we have when $m$ is sufficiently large that
\begin{align*}
    \begin{split}
        \frac{\Vert x_1^t-u^* \Vert^2}{\Vert x_1^0-u^* \Vert^2} \geq \frac{\sum_{k=t+1}^m q^{2k}}{\sum_{k=1}^m q^{2k}}=q^{2t}\frac{1-q^{2m-2t}}{1-q^{2m}} \geq \frac 12 q^{2t} \ ,
    \end{split}
\end{align*}
where the last inequality uses $m=2t$, $0<q<1$ and the fact that $q^{4t}-2q^{2t}+1=(q^{2t}-1)^2\geq 0$. 


Similar result holds for $x_2$:
\begin{equation}
    \frac{\Vert x_2^t-u^* \Vert^2}{\Vert x_1^0-u^* \Vert^2} \geq \frac 12 q^{2t} \ .
\end{equation}
To sum up, for any $t\in \mathbb N$ there exists an integer $m\geq 2t$ such that 
\begin{align}\label{eq:eq-lb}
    \begin{split}
        \text{dist}(z^t,\mathcal Z^*)^2 & \ \geq \text{dist}(x_1^t,u^*)^2+\text{dist}(x_2^t,u^*)^2 \\
        & \ \geq \frac{1}{2} q^{2t} \text{dist}(x_1^0,u^*)^2+\frac{1}{2} q^{2t} \text{dist}(x_2^0,u^*)^2\\
        & \ = \frac 12 q^{2t} \text{dist}\pran{x^0,(u^*,u^*)}^2\\
        & \ = \frac{1}{2}\pran{1-\frac{\alpha\sqrt 2}{L}}^{2t} \text{dist}(z^0,\mathcal Z^*)^2\ ,
    \end{split}
\end{align}
where second inequality is from Equation \eqref{eq:eq-lb}. Last equality uses $y_1^*=y_2^*=y_1^0=y_2^0=0$.
Taking square root we finally reach
\begin{equation}
    \text{dist}(z^t,\mathcal Z^*)\geq \frac{1}{\sqrt 2}\pran{1-\frac{\alpha\sqrt{2}}{L}}^{t} \text{dist}(z^0,\mathcal Z^*)\ .
\end{equation}

\end{proof}
 
\begin{rem}
    RsEGM with the stochastic gradient oracle $F_\xi$ given in the proof (see \eqref{eq:eq-lower-bound-oracle}) has upper bound equal to $\widetilde{\mathcal O}\pran{\frac{\Vert A \Vert_2}{\sigma^+_{min}(A)}\pran{\log{\frac{R_0}{\epsilon}}}^3}$. Compared with the above lower bound for unconstrained bilinear problem, $\Omega\pran{\frac{\Vert A \Vert_2}{\sigma^+_{min}(A)}\log{\frac{R_0}{\epsilon}}}$, we conclude that RsEGM is tight up to logarithmic factors.
\end{rem}

\section{Stochastic Oracles}\label{sec:app}

In this section, we propose four stochastic oracles and apply the main results to compute their total flop counts to obtain an approximate solution to the unconstrained bilinear problems and standard-form LP.


\subsection{Unconstrained Bilinear Problems}\label{sec:unconstrained-bilinear}

We consider unconstrained bilinear problems $$\min_{x\in \mathbb R^n}\max_{y\in \mathbb R^m} y^TAx+c^Tx-b^Ty \ .$$ Without loss of generality, we drop the linear terms $c^Tx-b^Ty$ because this can be achieved by shifting the origin. Now we consider the 
unconstrained bilinear problems of the form $$\min_{x\in \mathbb R^n}\max_{y\in \mathbb R^m} y^TAx \ ,$$ where $A\in \mathbb R^{m\times n}$. For comparison sake, we assume $\text{nnz}(A)\ge m+n$.


Converting to \eqref{eq:poi} and calculating $F(z)$, we have
\begin{equation*}
F(z)=F(x,y)=\left( \begin{array}{c} A^Ty\\-Ax \end{array}\right),\; g_1(x)=g_2(y)=0 \ .
\end{equation*}

First, we restate the total flop count of the optimal deterministic first-order primal dual method (it is achieved by restart primal-dual algorithms \cite{applegate2021faster}). In order to obtain $\epsilon$ accuracy solution, the number of primal-dual iterations is $\mathcal{O}\pran{\frac{\Vert A \Vert_2}{\alpha}\log{\frac{R_0}{\epsilon}}}$, thus the total flop counts is (by noticing the flop count of one primal-dual iteration is $O(\text{nnz}(A))$)
$$\mathcal{O}\pran{\text{nnz}(A)\frac{\Vert A \Vert_2}{\alpha}\log{\frac{R_0}{\epsilon}}} \ .$$

To ease the comparison of the flop counts, we state the following simple fact (\cite{palaniappan2016stochastic}):
{\footnotesize \begin{equation}\label{eq:eq-fact}
    \max{\{\Vert A_{i\cdot} \Vert_2,\Vert A_{\cdot j} \Vert_2\}}\leq \Vert A \Vert_2 \leq \Vert A \Vert_F \leq \sqrt{\max{\{m,n\}}}\cdot \max{\{\Vert A_{i\cdot} \Vert_2,\Vert A_{\cdot j} \Vert_2\}} \leq \sqrt{\max{\{m,n\}}}\cdot \Vert A\Vert _2
\end{equation}}

Next, we describe four different stochastic oracles that satisfies Assumption \ref{ass:oracal}. We will discuss their flop counts and how they improve deterministic restarted methods in different regimes. 

The first two oracles are based row/column sampling, i.e. the stochastic oracle is constructed using one row and one column. The stochastic oracle is given by the following
\begin{equation*}
    F_{\xi}(z)= \left( \begin{array}{c} \frac{1}{r_i}A_{i\cdot}y_i \\ -\frac{1}{c_j}A_{\cdot j}x_j \end{array} \right),\; \mathbb{P}(\xi=(i,j))=r_ic_j \ ,
\end{equation*}
where $c_j\ge 0$ for $j=1,..n$ and $r_i\ge 0$ for $i=1,...m$ are parameters of the sampling scheme with $\sum_{j=1}^n c_j = \sum_{i=1}^n r_i=1$.
Notice that
\begin{align*}
    \begin{split}
        \mathbb E\left[\Vert F_\xi(z) \Vert^2\right]& \ =\sum_{i=1}^m r_i \Vert \frac{1}{r_i}A_{i\cdot}y_i \Vert^2+\sum_{j=1}^n c_j \Vert -\frac{1}{c_j}A_{\cdot j}x_j \Vert^2\\
        & \ = \sum_{i=1}^m \frac{1}{r_i} y_i^2 \Vert A_{i\cdot} \Vert^2+\sum_{j=1}^n \frac{1}{c_j} x_j^2 \Vert A_{\cdot j} \Vert^2\\
        & \ \leq \max_{i,j}{\left\{\frac{\Vert A_{i\cdot} \Vert^2}{r_i},\frac{\Vert A_{\cdot j} \Vert^2}{c_j}\right\}} \Vert z \Vert^2 \ .
    \end{split}
\end{align*}
Thus the Lipschitz constant of stochastic oracle $F_\xi$ is bounded by
\begin{equation}\label{eq:eq-row-column}
    L\leq \sqrt{\max_{i,j}{\left\{\frac{\Vert A_{i\cdot} \Vert^2}{r_i},\frac{\Vert A_{\cdot j} \Vert^2}{c_j}\right\}}} \ .
\end{equation}

{\bf Oracle I: Uniform row-column sampling.}
In the first oracle, we uniformly randomly choose row and columns, i.e.,
\begin{equation*}
F_{\xi}(z)= \left( \begin{array}{c} \frac{1}{r_i}A_{i\cdot}y_i \\ -\frac{1}{c_j}A_{\cdot j}x_j \end{array} \right),\; \mathbb{P}(\xi=(i,j))=r_ic_j,\; r_i=\frac{1}{m},\; c_j=\frac{1}{n} \ .
\end{equation*}
Plugging in the choice of $r_i$ and $c_j$ into Equation \eqref{eq:eq-row-column}, the Lipschitz constant of $F_\xi$ is upper bounded by
\begin{equation}\label{eq:eq-l1}
    L_1\leq \sqrt{\max{\{ m\Vert A_{i\cdot} \Vert_2^2},n\Vert A_{\cdot j} \Vert_2^2 \}} \leq \sqrt{\max{\{m,n\}}}\cdot \max{\left\{\Vert A_{i\cdot} \Vert_2,\Vert A_{\cdot j} \Vert_2\right\}} \ .
\end{equation}

Combine Equation \eqref{eq:bd-oracle-complexity} together with \eqref{eq:eq-l1} and set $p=\frac{m+n}{\text{nnz}(A)}$. The total flop count of RsEGM to find an $\epsilon$-optimal solution using Oracle I is upper bounded by
{\small \begin{equation*} \widetilde{\mathcal O}\pran{\text{nnz}(A)\log{\frac{R_0}{\epsilon}}+\sqrt{\text{nnz}(A)(m+n)}\frac{\sqrt{\max{\{m,n\}}}\max{\{\Vert A_{i\cdot} \Vert_2,\Vert A_{\cdot j}\Vert_2}\}}{\alpha}\pran{\log{\frac{R_0}{\epsilon}}}^3} \end{equation*}}

Up to logarithmic factors, the RsEGM with uniformly sampled stochastic oracle is no worse than the deterministic method when A is dense according to Equation \eqref{eq:eq-fact}. And it improves the deterministic method when A is dense and $\max{\{\Vert A_{i\cdot} \Vert_2,\Vert A_{\cdot j}\Vert_2\}} \ll \Vert A \Vert_2$.

{\bf Oracle II: Importance row-column sampling.}
In the second oracle, we use importance sampling to set the probability to select rows and columns:
\begin{equation*}
F_{\xi}(z)= \left( \begin{array}{c} \frac{1}{r_i}A_{i\cdot}y_i \\ -\frac{1}{c_j}A_{\cdot j}x_j \end{array} \right),\; \mathbb{P}(\xi=(i,j))=r_ic_j,\; r_i=\frac{\Vert A_{i\cdot} \Vert_2^2}{\Vert A \Vert_F^2},\; c_j=\frac{\Vert A_{\cdot j} \Vert_2^2}{\Vert A \Vert_F^2} \ .
\end{equation*}
Plugging in the choice of $r_i$ and $c_j$ into Equation \eqref{eq:eq-row-column}, the Lipschitz constant of $F_\xi$ is upper bounded by
\begin{equation}\label{eq:eq-l2}
    L_2\leq\Vert A \Vert_F \ .
\end{equation}

Combining equation \eqref{eq:bd-oracle-complexity} with \eqref{eq:eq-l2} and setting $p=\frac{m+n}{\text{nnz}(A)}$, the total flop count of RsEGM with Oracle II to find an $\epsilon$-optimal solution is
$$\widetilde{\mathcal O}\pran{\text{nnz}(A)\log{\frac{R_0}{\epsilon}}+\sqrt{\text{nnz}(A)(m+n)}\frac{\Vert A \Vert_F}{\alpha}\pran{\log{\frac{R_0}{\epsilon}}}^3} \ .$$

It follows from \eqref{eq:eq-fact} that if $A$ is dense, RsEGM with importance sampled stochastic oracle is no worse than the optimal deterministic method up to a logarithmic factor. Furthermore, it improves the deterministic method when the stable rank $\frac{\Vert A \Vert_F^2}{\Vert A \Vert_2^2}\ll \frac{\text{nnz}(A)}{m+n}$, which is the case for low-rank dense matrix $A$.

In the above two oracles, the flop cost of computing the stochastic oracle is $\mathcal O (m+n)$. This can be further improved by using coordinate based sampling, where the flop cost per iteration is $\mathcal O(1)$ {(See Appendix \ref{sec:efficient} for details on efficient implementation of a more generic coordinate oracle)}.


{\bf Oracle III: Coordinate gradient estimator \cite{carmon2020coordinate}. }
The third oracle only updates a primal coordinate and a dual coordinate, thus the cost per iteration can be as low as $O(1)$. The idea follows from \cite{carmon2020coordinate}. More formally, we set
\begin{equation}\label{eq:oracle3}
F_{\xi}(z)= \left(  \frac{A_{i^xj^x}y_{i^x}}{p_{i^x j^x}}e_{j^x}, -\frac{A_{i^yj^y}x_{j^y}}{q_{i^y j^y}}e_{i^y}  \right)^T \ ,
\end{equation}
where $i_x, i_y\in\{1,...,n\}$, $j_x,j_y\in\{1,...,m\}$, and 
{\small \begin{equation*}
\mathbb{P}((i_x,j_x)=(i,j))=p_{ij}=\frac{\Vert A_{i\cdot} \Vert_1^2}{\sum_i \Vert A_{i\cdot} \Vert_1^2}\frac{|A_{ij}|}{\Vert A_{i\cdot} \Vert_1},\; \mathbb{P}((i_y,j_y)=(i,j))=q_{ij}=\frac{\Vert A_{\cdot j} \Vert_1^2}{\sum_j \Vert A_{\cdot j} \Vert_1^2}\frac{|A_{ij}|}{\Vert A_{\cdot j} \Vert_1}
\end{equation*}}
It is easy to check that
$\frac{A_{i^xj^x}y_{i^x}}{p_{i^x j^x}}e_{j^x}$ and $-\frac{A_{i^yj^y}x_{j^y}}{q_{i^y j^y}}e_{i^y}$ are unbiased estimators for $A^Ty$ and $-Ax$ respectively.
Furthermore, one can show the Lipschitz constant of the Oracle III (in expectation) is (see \cite{carmon2020coordinate} for more details):
\begin{equation}\label{eq:eq-l3}
    L_3=\max\left\{\sqrt{\sum_i \Vert A_{i\cdot} \Vert_1^2},\sqrt{\sum_j \Vert A_{\cdot j} \Vert_1^2}\right\}\leq \sqrt{m+n}\left\Vert \big|A\big| \right\Vert_2 \ ,
\end{equation}
where $\big|A\big|$ is the matrix with entry-wise absolute value of $A$. Thus Oracle III satisfies Assumption \ref{ass:oracal} with Lipschitz constant $L_3$.

Combine Equation \eqref{eq:bd-oracle-complexity} together with \eqref{eq:eq-l3} and set $p=\frac{1}{\text{nnz}(A)}$. The total cost of RsEGM with Oracle III equals
$$\widetilde{\mathcal O}\pran{\text{nnz}(A)\log{\frac{R_0}{\epsilon}}+\sqrt{\text{nnz}(A)}\frac{\max{\{\sqrt{\sum_i \Vert A_{i\cdot} \Vert_1^2},\sqrt{\sum_j \Vert A_{\cdot j} \Vert_1^2}}\}}{\alpha}\pran{\log{\frac{R_0}{\epsilon}}}^3}$$

For entry-wise non-negative dense matrix $A$,  RsEGM with Oracle III improve deterministic methods by a factor of $\sqrt{\max{\{m,n\}}}$. 
Unfortunately, there is no guaranteed improvement for a general matrix $A$ that involves negative value due to the definition of Lipschitz constant $L_3$.

{\bf Oracle IV: Coordinate gradient estimator with a new probability distribution.}
Inspired by Oracle III, we propose Oracle IV, where we utilize the same gradient estimation \eqref{eq:oracle3} but with different probability values $p_{ij}$ and $q_{ij}$: 
\begin{equation*}
\mathbb{P}((i_x,j_x)=(i,j))=p_{ij}=\frac{A_{ij}^2}{\Vert A \Vert _F^2},\; \mathbb{P}((i_y,j_y)=(i,j))=q_{ij}=\frac{A_{ij}^2}{\Vert A \Vert _F^2} \ .
\end{equation*}
Similarly, we know $\frac{A_{i^xj^x}y_{i^x}}{p_{i^x j^x}}e_{j^x}$ and $-\frac{A_{i^yj^y}x_{j^y}}{q_{i^y j^y}}e_{i^y}$ are unbiased estimators for $A^Ty$ and $-Ax$ respectively. Furthermore,
the Lipschitz constant of Oracle IV can be computed by:
\begin{align*}
    \begin{split}
        \mathbb E\left[\Vert F_\xi(z) \Vert^2\right]& \ =\sum_{i^x,j^x} p_{i^xj^x} (\frac{A_{i^xj^x}y_{i^x}}{p_{i^x j^x}})^2+\sum_{i^y,j^y} q_{i^y j^y} (-\frac{A_{i^yj^y}x_{j^y}}{q_{i^y j^y}})^2\\
        & \ =\sum_{i^x,j^x} \frac{(A_{i^xj^x})^2}{p_{i^xj^x}} (y_{i^x})^2+\sum_{i^y,j^y} \frac{(A_{i^yj^y})^2}{q_{i^y j^y}} (x_{j^y})^2\\
        & \ = \Vert A \Vert_F^2 \Vert z \Vert^2 \ .
    \end{split}
\end{align*}
Thus, $L_4=\Vert A \Vert_F$.
Combine Equation \eqref{eq:bd-oracle-complexity} with $L_4=\Vert A \Vert_F$ and set $p=\frac{1}{\text{nnz}(A)}$. The total flop cost of RsEGM with Oracle IV to find an $\epsilon$-solution  becomes $$\widetilde{\mathcal O}\pran{\text{nnz}(A)\log{\frac{R_0}{\epsilon}}+\sqrt{\text{nnz}(A)}\frac{\Vert A \Vert_F}{\alpha}\pran{\log{\frac{R_0}{\epsilon}}}^3} \ .$$
Note that when the matrix $A$ in the unconstrained bilinear problem is dense, RsEGM improves the total flop counts of deterministic primal-dual method (upto log terms) by at least a factor of $\sqrt{\max{\{m,n\}}}$. The improvement does NOT require extra assumptions on the spectral property of $A$ (such as Oracle II) nor non-negativity of entries of $A$ (such as Oracle III).

\subsection{Linear programming}\label{sec:oracle-lp}
Consider the primal-dual formulation of standard form LP \eqref{eq:lp}. For the ease of comparison, we assume $\text{nnz}(A)\ge m+n$. Furthermore, we set
\begin{equation}
F(z)=\left( \begin{array}{c} A^Ty\\-Ax \end{array}\right),\; g_1(x)=c^Tx+\iota_{\{x\geq 0\}},\; g_2(y)=b^Ty.
\end{equation}
Recall that $R_0=\text{dist}(z^{0,0},\mathcal{Z}^*)$. The total flop count of restarted primal-dual algorithms is (see \cite{applegate2021faster} for details)
$$\mathcal{O}\pran{\text{nnz}(A)\frac{\Vert A \Vert_2}{1/[H(1+\Vert z^{0,0} \Vert+R_0)]}\log{\frac{R_0}{\epsilon}}}$$

\textbf{(i)} Consider RsEGM with stochastic Oracle II (importance sampling)
\begin{equation*}
F_{\xi}(z)= \left( \begin{array}{c} \frac{1}{r_i}A_{i\cdot}y_i \\ -\frac{1}{c_j}A_{\cdot j}x_j \end{array} \right),\; \mathbb{P}(\xi=(i,j))=r_ic_j,\; r_i=\frac{\Vert A_{i\cdot} \Vert_2^2}{\Vert A \Vert_F^2},\; c_j=\frac{\Vert A_{\cdot j} \Vert_2^2}{\Vert A \Vert_F^2} \ .
\end{equation*}
Similar to the unconstrained bilinear case (Section \ref{sec:unconstrained-bilinear}), we can obtain the flop count 
$$\widetilde{\mathcal O}\pran{\text{nnz}(A)\log{\frac{R_0}{\epsilon}}+\sqrt{\text{nnz}(A)(m+n)}\frac{\Vert A \Vert_F}{1/[H(1+\Vert z^{0,0} \Vert+R_0)]}\pran{\log{\frac{R_0}{\epsilon}}}^5}$$

\textbf{(ii)} Consider RsEGM with Oracle IV, i.e.,
\begin{equation*}
F_{\xi}(z)= \left(  \frac{A_{i^xj^x}y_{i^x}}{p_{i^x j^x}}e_{j^x}, -\frac{A_{i^yj^y}x_{j^y}}{q_{i^y j^y}}e_{i^y}  \right)^T,\; p_{ij}=\frac{A_{ij}^2}{\Vert A \Vert _F^2},\; q_{ij}=\frac{A_{ij}^2}{\Vert A \Vert _F^2} \ .
\end{equation*}

{
With a similar calculation (Section \ref{sec:unconstrained-bilinear}), we can obtain the flop count:
$$\widetilde{\mathcal O}\pran{\text{nnz}(A)\log{\frac{R_0}{\epsilon}}+\sqrt{\text{nnz}(A)}\frac{\Vert A \Vert_F}{1/[H(1+\Vert z^{0,0} \Vert+R_0)]}\pran{\log{\frac{R_0}{\epsilon}}}^5} \ .$$
Note that $\Vert A \Vert_F\le \sqrt{\text{nnz}(A)}\|A\|_2$, thus RsEGM always has a complexity no worse than the deterministic primal-dual methods (upto log terms). In the case when the constraint matrix $A$ is dense, the improvement of the total flop counts over its deterministic counterpart is by at least a factor of $\sqrt{\max{\{m,n\}}}$.
}


\section{{Numerical Experiment}}

In this section, we present numerical experiments on two classes of problems: matrix games and linear programming.

\subsection{Matrix games}\label{sec:mg-numerical}
{\bf Problem instances.} Matrix games can be formulated as a primal-dual optimization problem with the following form:
\begin{equation*}
    \min_{x\in \Delta_n}\max_{y\in \Delta_m} y^TAx \ ,
\end{equation*}
where $\Delta_n$ and $\Delta_m$ are $n$-dimensional and $m$-dimensional simplex respectively and matrix $A\in \mathbb R^{m\times n}$. In the experiment, we generate three matrix game instances with scale $m=n=1000$ using the instance-generation code from \cite{alacaoglu2021stochastic}, which generated matrix game instances from classic literature \cite{nemirovski2009robust,nemirovski2013mini}. 

{\bf Progress metric.} We use duality gap, i.e., $\max_{\hat x\in \Delta_n,\hat y\in \Delta_m} \hat y^TAx-y^TA\hat x$ as the performance measure for different algorithms. The duality gap can be easily evaluated by noticing it equals $\max_i(Ax)_i-\min_j(A^Ty)_j$. 

{\bf Projection and stochastic oracle.} We use the projection oracle developed in \cite{condat2016fast} for the projection step onto standard simplex. We utilize the importance row-column sampling (i.e., Oracle II in Section \ref{sec:unconstrained-bilinear}) to implement RsEGM. 

\begin{figure}[h!]
	\centering
	\begin{tabular}{l c c c c}
		\hspace{-2Cm}
		& \includegraphics[width=0.33\textwidth]{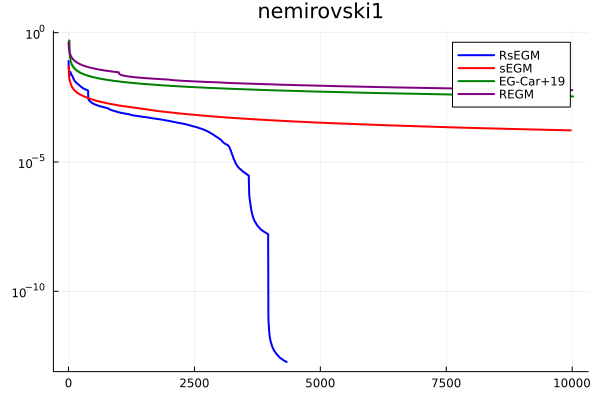}
		\hspace{-0.5Cm}
		& \includegraphics[width=0.33\textwidth]{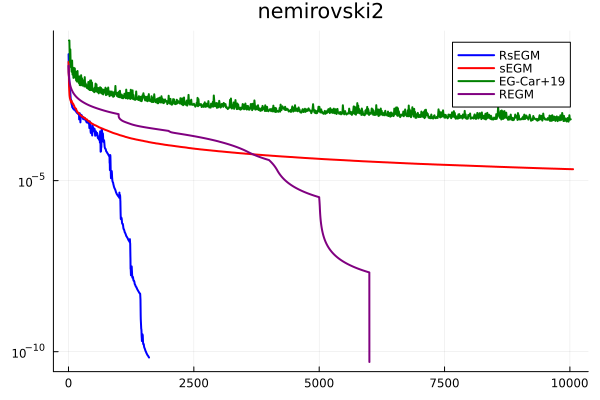}
        \hspace{-0.5Cm}
        & \includegraphics[width=0.33\textwidth]{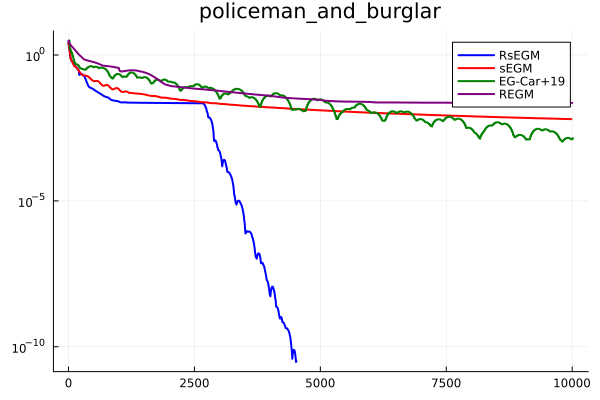}
	\end{tabular}
	\caption{Convergence behaviour on matrix games with simplex constraint.}
	\label{fig:mg}
\end{figure}

{\bf Results.} We compare RsEGM with several notable methods in literature for solving matrix games with simplex constraint, namely, stochastic extragradient method with variance reduction (sEGM) \cite{alacaoglu2021stochastic}, variance-reduced method for matrix games (EG-Car+19) \cite{carmon2019variance} and deterministic restarted extragradient method (REGM) \cite{applegate2021faster}.
Figure \ref{fig:mg} presents the numerical performance of different methods on the three matrix game instances. The $x$-axis is the number of iterations for a deterministic gradient step (a stochastic gradient step is adjusted proportionally), and the $y$-axis is the duality gap. We can observe that the proposed RsEGM exhibits the fastest convergence to obtain a high-accuracy solution and clearly exhibits a linear convergence rate, which verifies our theoretical guarantees. Furthermore, we can observe that both sEGM and EG-Car+19 exhibit slow convergence compared to RsEGM. Moreover, although REGM also enjoys linear convergence theoretically, the linear rate is only shown up in instance nemirovski2 and has not yet been triggered in the other two instances within 10000 iterations.  

\subsection{Linear programming}
{\bf Problem instances.} In the experiments, we utilize the instances from \href{https://miplib.zib.de/tag_collection.html}{\texttt{MIPLIB 2017}} to compare the performance of different algorithms. We randomly select six instances from \href{https://miplib.zib.de/tag_collection.html}{\texttt{MIPLIB 2017}} and look at its root-node LP relaxation. For these instances, we first convert the instances to the following form:
\begin{align}\label{eq:practical-lp-form}
    \begin{split}
        & \min_x\; c^Tx\\
        & \ \mathrm{s.t.}\; A_{E}x=b_E \\
        & \ \quad\;\; A_{I}x\leq b_I\\
        & \ \quad\;\; x\geq 0 \ ,
    \end{split}
\end{align}
and consider its primal-dual formulation
\begin{align}\label{eq:practical-lp}
    \begin{split}
        \min_{x\geq 0}\max_{y_I\leq 0}\; c^Tx-y^TAx+b^Ty \ , 
    \end{split}
\end{align}
where $A_E\in\mathbb R^{m_E\times n}$, $A_I\in \mathbb R^{m_I\times n}$, $A=\begin{pmatrix}
    A_{E} \\ A_{I}
\end{pmatrix}\in\mathbb R^{(m_E+m_I)\times n}$ and $b=\begin{pmatrix}
    b_E \\ b_I
\end{pmatrix}\in\mathbb R^{(m_E+m_I)}$. We then compare the numerical performance of different primal-dual algorithms for solving \eqref{eq:practical-lp}. This is the problem format supported by the PDHG-based LP solver PDLP~\cite{applegate2021practical}.

{\bf Progress metric.} We use KKT residual of \eqref{eq:practical-lp}, i.e., a combination of primal infeasibility, dual infeasibility and primal-dual gap, to measure the performance of current iterates.  More formally, the KKT residual of \eqref{eq:practical-lp-form} is given by 
\begin{equation*}
    \mathrm{KKT}(x,y)=\left\Vert \begin{pmatrix}
        A_{E}x-b_E \\ [A_{I}x-b_{I}]^+ \\ [-x]^+  \\ [A^Ty-c]^+ \\ [y_I]^+ \\  [c^Tx-b^Ty]^+
    \end{pmatrix}  \right\Vert_2 \ .
\end{equation*}

{\bf Stochastic oracles.} Similar to matrix games, we utilize the importance row-column sampling scheme to implement RsEGM, namely, stochastic Oracle II presented in Section \ref{sec:oracle-lp}.

{\bf Results.} We compare RsEGM with two methods for solving linear programming: stochastic extragradient method with variance reduction without restart (sEGM) \cite{alacaoglu2021stochastic} and deterministic restarted extragradient method (REGM) \cite{applegate2021faster}.  Note that the variance-reduced method for matrix games (EG-Car+19) \cite{carmon2019variance} compared in Section \ref{sec:mg-numerical} is not considered for linear programming since it requires bounded domain, which is in general not satisfied by real-world LP instances, as those in the \texttt{MIPLIB 2017} dataset.

\begin{figure}[ht!]
	\centering
	\begin{tabular}{l c c c c}
		\hspace{-2Cm}
		& \includegraphics[width=0.33\textwidth]{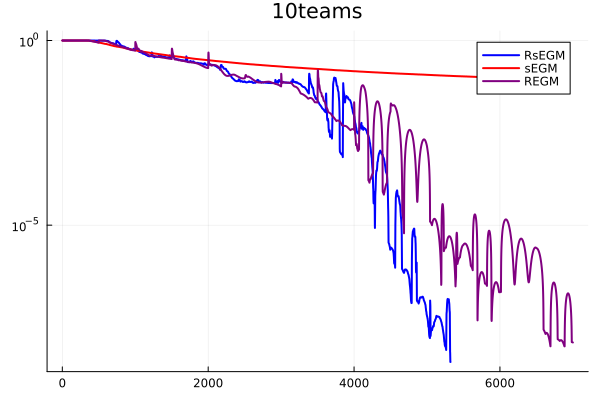}
		\hspace{-0.5Cm}
		& \includegraphics[width=0.33\textwidth]{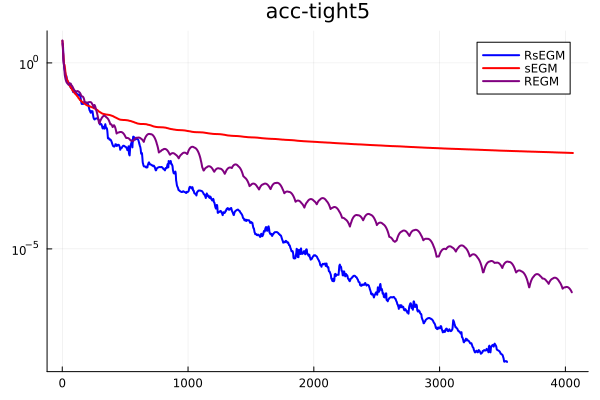}
        \hspace{-0.5Cm}
        & \includegraphics[width=0.33\textwidth]{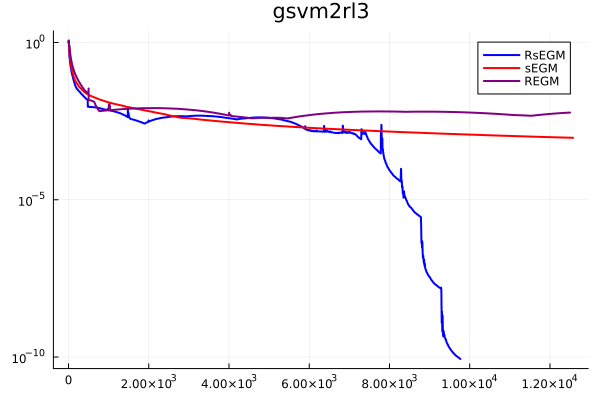}\\
        \hspace{-2Cm}
        & \includegraphics[width=0.33\textwidth]{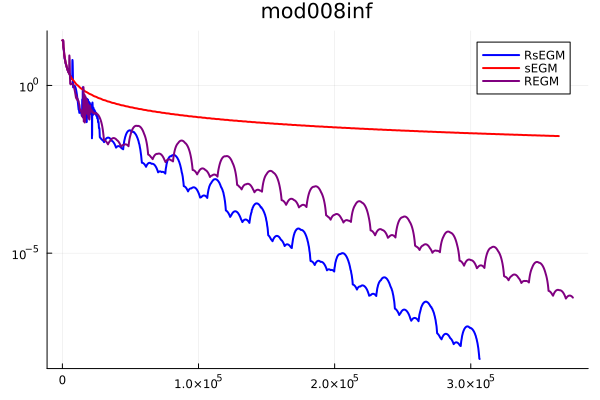}
		\hspace{-0.5Cm}
		& \includegraphics[width=0.33\textwidth]{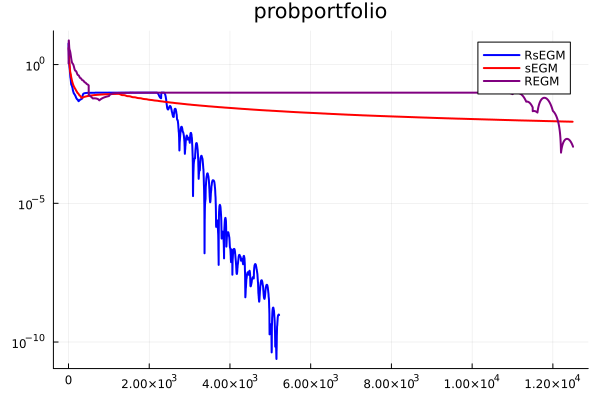}
        \hspace{-0.5Cm}
        & \includegraphics[width=0.33\textwidth]{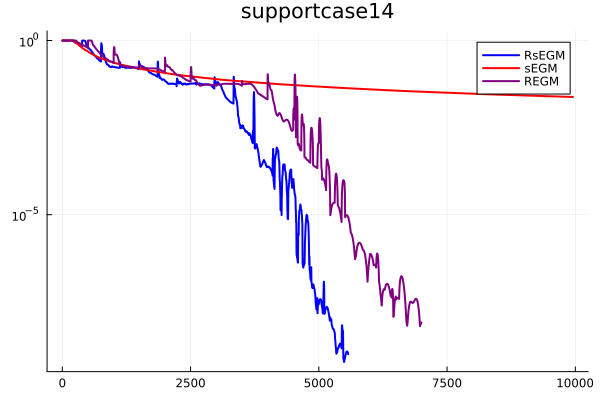}\\
		\hspace{-1Cm}
	\end{tabular}
	\caption{Convergence behaviour on LP instances from \texttt{MIPLIB}.}
	\label{fig:lp}
\end{figure}

Figure \ref{fig:lp} presents the convergence behaviors of different methods on the six instances. The $x$-axis is the number of iterations for a deterministic gradient step (a stochastic gradient step is adjusted proportionally), and the $y$-axis is the KKT residual. Again, we can see that RsEGM exhibits the fastest convergence to achieve the desired accuracy $10^{-5}$, and sEGM has slow sublinear convergence. REGM has competitive performance with linear convergence, but RsEGM enables earlier triggers of linear convergence and a faster eventual rate in general.

\section{Conclusions}
In this work, we introduce a stochastic algorithm for sharp primal-dual problems, such as linear programming and bilinear games, using variance reduction and restart. We show that our proposed stochastic methods enjoy a linear convergence rate, which improves the complexity of the existing algorithms in the literature. {We also show that our proposed algorithm achieves the optimal convergence rate (upto a log term) in a wide class of stochastic algorithms.}

\section*{{Acknowledgement}}
The authors would like to thank Warren Schudy for proposing and discussing the efficient update of the coordinate stochastic oracle (Oracle IV) for LP. 

\bibliographystyle{amsplain}
\bibliography{ref-papers}

\newpage

\appendix

\section{Appendix}

\subsection{Sharpness and metric sub-regularity}\label{sec:sharp-subreg}
In this subsection we establish a connection between sharpness based on the normalized duality gap \cite{applegate2021faster} and metric sub-regularity \cite{latafat2019new}. The result is summarized in the following proposition.

\begin{prop}\label{thm:thm-sharp-subreg}
    Consider primal-dual problem $\min_{x\in \mathcal{X}}\max_{y\in \mathcal{Y}}\mathcal L(x,y)$. If for any $z\in \mathcal{U}(z^*)$ the primal-dual problem is sharp with constant $\alpha$, i.e. $\alpha \text{dist}(z, \mathcal Z^*)\leq \rho_r(z)$ for all $r>0$, then it satisfies metric sub-regularity at $z^*$ for 0, i.e. $\alpha \text{dist}(z,\mathcal Z^*)\leq \text{dist}(0,\partial \mathcal L(z))=\inf_{(u,-v)\in \partial \mathcal L(z)}\Vert (u,-v) \Vert,\; \forall z\in \mathcal{U}(z^*)$.
\end{prop}

{\bf Proof.}
For any $(u,-v)\in \partial \mathcal L(z)$
\begin{align*}
    \begin{split}
        \alpha \text{dist}(z, \mathcal Z^*) & \leq \rho_r(z)= \frac 1r \max_{\hat{z}\in W_r(z)}{\mathcal L(x,\hat y)-\mathcal L(\hat x,y)}\\
        & \leq \frac 1r \max_{\hat{z}\in B_r(z)}{\mathcal L(x,\hat y)-\mathcal L(\hat x,y)}\\
        & = \frac 1r \max_{\hat{z}\in B_r(z)}{\mathcal L(x,\hat y)-\mathcal L(x,y)+\mathcal L(x,y)-\mathcal L(\hat x,y)}\\
        & \leq \frac 1r \max_{\hat{z}\in B_r(z)}{v^T(\hat y -y)+(-u)^T(\hat x -x)}\\
        & = \frac 1r \max_{\hat{z}\in B_r(z)}{-(u,-v)^T(\hat z -z)}\\
        & = \frac 1r r\Vert (u,-v) \Vert = \Vert (u,-v) \Vert \ , 
    \end{split}
\end{align*}
where the first inequality utilizes the definition of sharpness, the third one is due to the convexity-concavity of $\mathcal L(x,y)$.

Take infimum over $\partial \mathcal L(z)$ we have
\begin{equation*}
    \alpha \text{dist}(z, \mathcal Z^*) \leq  \inf_{(u,-v)\in \partial \mathcal L(z)}\Vert (u,-v) \Vert=\text{dist}(0,\partial \mathcal L(z))
\end{equation*}
\qed

\subsection{LP does not satisfy metric sub-regularity globally}\label{sec:lp_not_global_sharp}

In this subsection, we present a counter-example to show that the Lagrangian's generalized gradient of LP does not satisfy metric sub-regularity globally. 

Consider an LP
\begin{equation*}
    \min_{x\in \mathbb R^2}{\mathbf{1}^Tx},\quad \text{s.t.}\; Ax\leq b, x\geq 0
\end{equation*}
where $A=(0\;1)\in \mathbb R^{1\times 2}$ and $b=1$. The optimal solution is unique $x^*=(0,0)$. The primal dual form is
\begin{equation*}
    \mathcal L(x,y)=\mathbf{1}^Tx+y^T(Ax-b)+\iota_{\{x\geq 0\}}-\iota_{\{y\geq 0\}}\ ,
\end{equation*}
and the optimal primal-dual solution $z^*=(\mathbf 0_2,0)$. Notice that
\begin{equation*}
    \partial L(z)=\bigg\{(p,q)\in\mathbb R^{2+1}\bigg|p\in\mathbf{1} +A^Ty+\partial \iota_{\{x\geq 0\}},q\in-Ax+b+\partial\iota_{\{y\geq 0\}}\bigg\}
\end{equation*}
thus
\begin{equation*}
    \text{dist}(0,\partial \mathcal L(z))^2 = \inf_{(p,q)\in\partial L(z)}{\Vert p \Vert^2+q^2}\leq 1+(1+y)^2+(1-x_2)^2 
\end{equation*}

If the Lagrangian's generalized gradient of LP satisfies global metric sub-regularity, i.e., there exists a constant $\eta$ such that $\eta\text{dist}(z,\mathcal Z^*)\leq \text{dist}(0,\partial \mathcal L(z)), \; \forall z\in \mathcal Z$, then
\begin{equation*}
    x_1^2+x_2^2+y^2\leq \frac{1}{\eta^2} (1+(1+y)^2+(1-x_2)^2),\; \forall (x_1,x_2,y)\in\mathbb R^3_+\ .
\end{equation*}
However, this is impossible by fixing $x_2, y$ and letting $x_1\rightarrow \infty$. This leads to a contradiction.

Therefore LP does not satisfy the global metric sub-regularity. Although we use Euclidean norm for illustration, due to the equivalence of norms on $\mathbb R^n$, LP does not satisfy the global metric sub-regularity under any norm. Notice that metric subregularity is a weaker condition than sharpness (see Appendix \ref{sec:sharp-subreg}). A direct consequence is that LP is not globally sharp.




\subsection{{Efficient update for coordinate sampling}}\label{sec:efficient}
In Section \ref{sec:unconstrained-bilinear}, we propose two coordinate gradient estimators, Oracle III and Oracle IV. Here, we present an efficient computation approach that $\mathcal{O}(\text{nnz}(A))$ flops for a snapshot step (i.e., the steps when updating the snapshot $w$), and $\mathcal{O} (1)$ flops for a non-snapshot step (i.e., the steps between two snapshot steps). With such efficient update rules, we obtain the total cost computation for the Oracle IV rows in Table \ref{tab:bilinear} and Table \ref{tab:lp}.


Consider a generic coordinate stochastic oracle $$F_{\xi}(z)=(c_xe_x,d_ye_y) \ ,$$ where $e_x,e_y$ are standard unit vectors in $\RR^m$ and $\RR^n$, respectively, and $c_x, d_y$ are two arbitrary scalars. Oracle III and IV are special cases for this generic coordinate stochastic oracle.

Let $w^{k^j}$ be the $j$-th snapshot step, i.e., the $j$-th updates of $w$. Consider a non-snapshot step $k$ with $k^j\le k< k^{j+1}-1$, i.e., the steps between two consecutive snapshot steps, then we know $w^k\equiv w^{k^j}$ and $u=p(w^{k^j})-\tau F(w^{k^j})$ do not change for any $k^j<k< k^{j+1}-1$. Furthermore, it is easy to check that $F_{\xi_{k}}(z^{k+1/2})_i=F_{\xi_{k}}(w^{k})_i$ has update
\begin{equation*}
    z^{k+1}_{i}=\begin{cases}
        \max\{0,u_i/p+(1-p)(z^{k}_{i}-u_i/p)\}, & \mathrm{if\ coordinate}\ i \mathrm{\ is \ constrained\ with\ }z_i\geq 0\\
        u_i/p+(1-p)(z^{k}_{i}-u_i/p), & \mathrm{if\ coordinate}\ i \mathrm{\ is \ unconstrained}
    \end{cases} \ .
\end{equation*}

Generalizing this simple observation, a closed-form update rule holds for iterates $z^{k'}$ with $k^j\le k < k' \le k^{j+1}$ that
\begin{equation}\label{eq:update-coordinate}
    z^{k'}_{i}=\begin{cases}
        \max\{0,u_i/p+(1-p)^{k'-k}(z^{k}_{i}-u_i/p)\}, & \mathrm{if\ coordinate}\ i \mathrm{\ is \ constrained\ with\ }z_i\geq 0\\
        u_i/p+(1-p)^{k'-k}(z^{k}_{i}-u_i/p), & \mathrm{if\ coordinate}\ i \mathrm{\ is \ unconstrained}
    \end{cases} \ .
\end{equation}

Thus motivated, we can run Algorithm \ref{alg:segm} with an efficient update rule using the following procedure:


\textbf{Snapshot steps.} In a snapshot step (i.e., $k=k^j$), we compute $z^{k^j}$ using \eqref{eq:update-coordinate}, set the iteration of last update of coordinate $i$ (denoted as $k_i$) to $k^j$, namely, $k_i=k^j$ for $1\le i\le m+n$, record $z^{k^j}$ in memory, and set $w^{k^j}=z^{k^j}$. The computational cost is $\mathcal O(m+n)$. Then we compute $u=pw^{k^j}-\tau F(w^{k^j})$ and store $u$ in memory. The computational cost is $\mathcal O(\text{nnz}(A))$. Finally, we compute the running average of iterates $\widetilde{z}^{k^j}$, and the cost is $\mathcal O(m+n)$ by taking advantage of the sum over geometric series.



\textbf{Non-snapshot steps.} In a non-snapshot step (i.e, $k^j<k<k^{j+1}-1$), if $F_{\xi_k}(z_{k+1/2})$ updates coordinate $i$ and $F_{\xi_k}(z^{k+1/2})_i\neq F_{\xi_k}(w^{k})_i$, we compute $z^{k}_{i}$ using \eqref{eq:update-coordinate} and update $z^{k+1/2}_{i}$ and $z^{k+1}_{i}$. Then set $k_i=k+1$. We store $z^{k_i}_{i}$ in memory. The cost is $\mathcal O(1)$.


As such, the cost when updating the snapshot $w$ is $\mathcal{O}(\text{nnz}(A))$ and the cost for a non-snapshot step is $\mathcal{O}(1)$ for Oracle III and Oracle IV discussed in Section \ref{sec:app}.

\subsection{Comparison with SPDHG}\label{sec:compare-spdhg}
\cite[Theorem 4.6]{alacaoglu2019convergence} shows the linear convergence of SPDHG for problems with global metric sub-regularity, such as the unconstrained bilinear problem. However, the obtained linear convergence rate may involve parameters that are not easy to interpret. We here apply their results to unconstrained bilinear problems with more transparent parameters so that we make a comparison with our Theorem \ref{thm:thm-global}. 

Using the notation in \cite{alacaoglu2019convergence}, they obtained
\begin{equation*}
    \mathbb E\left [ \frac{C_1}{2}\Vert x^k-x_*^k \Vert_{\tau^{-1}}^2 +\frac 12 \Vert y^{k+1}-y_*^{k+1} \Vert_{D(\sigma)^{-1}P^{-1}}^2 \right] \leq (1-\rho)^k2\Phi^0 \ ,
\end{equation*}
where $\rho=\frac{C_1\underline p}{2\zeta}$, $\underline p=\min_ip_i$, $C_1=1-\gamma$, $\zeta=2+\frac 2 {\alpha^2}\Vert H-N \Vert^2$. For simplicity, we consider the uniform sampling scheme with $p_i=\frac 1n$, and thus $\underline p=\min_ip_i=\frac 1n$ (the other schemes follow with a similar arguments).
The total number of iterations for SPDHG to achieve $\epsilon$-accuracy is of order
\begin{align*}
    \begin{split}
        & \ \mathcal O\pran{\frac 1\rho \log{\frac 1 \epsilon}}=  \ \mathcal O \pran{\frac{\zeta}{\underline p}\log \frac 1 \epsilon}
        = \ \mathcal O(n\zeta \log \frac 1 \epsilon)
        =  \ \mathcal O(n\frac{\Vert H-N \Vert^2}{\alpha^2}\log \frac 1 \epsilon)\ .
    \end{split}
\end{align*}
where the definitions of $H$ and $N$ are defined in \cite[Lemma 4.4]{alacaoglu2019convergence}.
The operators $H$ and $N$ are indeed a bit complicated so that the upper bound is hard to derive. Now we lower bound the term $\Vert H-N \Vert^2$. 
\begin{align*}
    \begin{split}
        \Vert H-N \Vert^2 & \ \geq \frac{\Vert (H-N)q \Vert_2^2}{\Vert q \Vert_2^2} \ \geq \frac{\Vert \tau^{-1}x(1) \Vert_2^2+\Vert Ax(1) \Vert_2^2}{\Vert x(1) \Vert_2^2}  \ \geq \frac{\Vert Ax(1) \Vert_2^2}{\Vert x(1) \Vert_2^2} =\Vert A \Vert^2 \ ,
    \end{split}
\end{align*}
where the first inequality is due to the definition of operator norm, the second inequality follows from the choice of $q=(x(1),0,...,0)$ and the last equality is from the choice of $x(1)$ such that $x(1)=\text{argmax}_u \frac{\Vert Au \Vert_2^2}{\Vert u \Vert_2^2}$.

Thus we have the complexity of SPDHG is at least
{\small \begin{equation*}
    \mathcal O\pran{\frac 1\rho \log{\frac 1 \epsilon}}=\mathcal O(n\frac{\Vert H-N \Vert^2}{\alpha^2}\log \frac 1 \epsilon) \geq \mathcal O(n\frac{\Vert A \Vert^2}{\alpha^2}\log \frac 1 \epsilon) \geq \mathcal O(\frac{\Vert A \Vert_F^2}{\alpha^2}\log \frac 1 \epsilon) =\mathcal O(\kappa_F^2\log \frac 1 \epsilon) \ .
\end{equation*}}
In contrast, the complexity of RsEGM is $\widetilde{\mathcal O}(\kappa_F\pran{\log \frac 1 \epsilon}^3)$. This showcase RsEGM has an improved linear convergence rate for solving unconstrained bilinear problems compared to SPDHG in terms of the condition number (after ignoring the $\log$ terms). Such an improvement comes from the restarted scheme, similar to the deterministic case~\cite{applegate2021faster}. Lastly, we would like to mention that SPDHG is a semi-stochastic algorithm, where the primal is a full gradient update and the dual is a stochastic gradient update, while RsEGM  takes stochastic gradient steps in both primal and dual space (except the snapshot steps).

\end{document}